\documentclass{siamart0516}
\usepackage{lipsum}
\usepackage{amssymb}
\usepackage{amsmath}
\usepackage{amsfonts}
\usepackage{graphicx}
\usepackage{epstopdf}
\usepackage{algorithmic}
\ifpdf
  \DeclareGraphicsExtensions{.eps,.pdf,.png,.jpg}
\else
  \DeclareGraphicsExtensions{.eps}
\fi
\numberwithin{equation}{section}

\topmargin by-\headsep  \textheight   236mm   
\def\to{\rightarrow}

\newcommand{\q}{\quad}   \newcommand{\qq}{\qquad}

 \def\bs{\bigskip} 

\def\beqn{\begin{eqnarray}}  \def\eqn{\end{eqnarray}}
\def\beqnx{\begin{eqnarray*}} \def\eqnx{\end{eqnarray*}}
\def\bn{\begin{enumerate}} \def\en{\end{enumerate}}

\newcommand{\TheTitle}{Mittag-Leffler stability of F-LMMs for F-ODEs} 
\newcommand{\TheAuthors}{Wang Dongling and Zou Jun}
\headers{\TheTitle}{\TheAuthors}

\title{Mittag--Leffler stability of numerical solutions \\ to time fractional ODEs
}

\author{
Dongling Wang\thanks{School of Mathematics and Computational Science, Xiangtan University, Xiangtan, Hunan 411105, P.R. China.
The work of this author was partially supported by National Natural Science Foundation of China (Grant No. 11871057, 91630205).
(\email{wdymath@xtu.edu.cn}).
}
  \and
Jun Zou\thanks{Department of Mathematics, The Chinese University of Hong Kong Shatin, N.T., Hong Kong.
 The work of this author was substantially supported by Hong Kong RGC General Research Fund (project 14306718). 
  (\email{zou@math.cuhk.edu.hk}).}
  }
 
\ifpdf
\hypersetup{
  pdftitle={\TheTitle},
  pdfauthor={\TheAuthors}
}
\fi

\begin{document}
\maketitle

\begin{abstract}
The asymptotic stable region and long-time decay rate of solutions to linear homogeneous Caputo time fractional ordinary 
differential equations (F-ODEs) are known to be completely determined by the eigenvalues of the coefficient matrix.
Very different from the exponential decay of solutions to classical ODEs, solutions of F-ODEs decay only polynomially, leading to the so-called Mittag-Leffler stability, which was already extended to semi-linear F-ODEs with small perturbations.   
This work is mainly devoted to the qualitative analysis of the long-time behavior of numerical solutions.
By applying the singularity analysis of generating functions developed by Flajolet and Odlyzko (SIAM J.\,Disc.\,Math.\,3\,(1990), 216-240), we are able to prove that both $\mathcal{L}$1 scheme and strong $A$-stable fractional linear multistep methods (F-LMMs) can preserve the numerical Mittag-Leffler stability for linear homogeneous F-ODEs exactly as in the continuous case. Through an improved estimate of the discrete fractional resolvent operator,
we show that strong $A$-stable F-LMMs are also Mittag-Leffler stable for semi-linear F-ODEs under small perturbations.
For the numerical schemes based on $\alpha$-difference approximation to Caputo derivative, 
we establish the Mittag-Leffler stability for semi-linear problems by making use of properties of the Poisson transformation
The new results and analyses provide not only the rigorous justifications and explanations of 
the Mittag-Leffler stability of numerical solutions with exact decay rate, but also establish 
some close connection between the continuous and discrete F-ODEs. 
Numerical experiments are presented for several typical time fractional evolutional equations, 
including time fractional sub-diffusion equations and semi-linear F-ODEs.
All the numerical results exhibit the typical long-time polynomial decay rate, which is fully consistent with our theoretical predictions.
\end{abstract}

 
\qq\qq \qq $\alpha$-difference.


\medskip
\section{Introduction}\label{Introduction}
Stability is one of the most fundamental issues for all time-dependent differential equations, 
and a deep understanding of the stability of linear problems is often 
a key step to the understanding of nonlinear models.
In this work, we are mainly concerned with the numerical stability of time fractional ordinary differential equations (F-ODEs) 
of the form (with $0<\alpha<1$): 
\begin{equation} \label{model}
\begin{split}
\mathcal{D}_{t}^{\alpha}y(t)=Ay+f(t, y), \quad t>0 
\end{split}
\end{equation}
for $y\in\mathbb{R}^{d}$ satisfying the initial value $y(0)=y_{0}$, where 
$A\in\mathbb{R}^{d\times d}$ is a real matrix, 
$f: {\mathbb{R}\times} \mathbb{R}^{d}\to \mathbb{R}^{d}$ is continuous,
and $\mathcal{D}_{t}^{\alpha}$ stands for the Caputo fractional derivative of order $\alpha$.  
Without loss of generality, we assume $f(t,0)=0$ so that the trivial solution $y=0$ 
is always an equilibrium solution to \eqref{model}.
The main difficulty in stability analysis of F-ODEs 
lies in the nonlocal nature of fractional derivatives. 
Let us first recall some definitions of stability for the trivial solution 
to the model \eqref{model} \cite{ Cong18,Cong19,Diethelm10}.

\begin{definition} \label{def:Cong}
The trivial solution of  F-ODEs \eqref{model} is said to be stable if for any $\varepsilon>0$, there exists $\delta=\delta(\varepsilon)>0$ such that for any
$\|y_{0}\|<\delta$ we have $\|y(t)\|\leq\varepsilon$ for all $t\geq 0$; and the trivial solution is said to be  asymptotically stable if it is stable and $\lim_{t\to\infty}\|y(t)\|=0$.  
It is said to be Mittag-Leffler stable if there exist positive constants $\beta, \delta$ and $M$ independent of $t$ such that 
\begin{equation} \label{MLstab}
\begin{split}
\sup_{t\geq 0} t^{\beta}\|y(t)\|\leq M~~\text{for any}~ \|y_{0}\|\leq \delta.
\end{split}
\end{equation}
\end{definition}
Here and in the sequel, we use $\|\cdot\|$ for the standard Euclidean norm
in $\mathbb{R}^{d}$.
It is known from \cite{Cong19} that the index $\beta$ in \eqref{MLstab} stays in the range 
$0<\beta\leq\alpha$ for the model equation \eqref{model} (cf.\,\cref{lem:Cong}).
There is an alternative definition of the Mittag-Leffler stability, by replacing 
the inequality \eqref{MLstab} by
$
\|y(t)\|\leq  V(y_{0}) E_{\alpha}(-Lt^{\alpha}),
$
where $L> 0$ and the function $V(y)$ is locally Lipschitz continuous and satisfies that $V(0)=0$ and $V(y)\geq0$. 
This alternative definition highlights the boundedness of the solutions by the Mittag-Leffler function 
$E_{\alpha}(z)$, which is similar to the exponential function in the classical ODEs. 

Basic definitions and properties of fractional calculus and Mittag-Leffler functions are included in \cref{sec:appendixA}.  
By means of the asymptotic expansion of Mittag-Leffler functions, 
the above two definitions are known to be essentially equivalent. 
Mittag-Leffler stability not only implies the asymptotic stability of the trivial solution to \eqref{model}, 
but also characterizes its long-time polynomial decay rate, which is 
an important common property of the solutions to F-ODEs.
Just like homogeneous linear ODEs with constant coefficients, the stability of the solutions to 
linear homogenous F-ODEs is completely determined by the eigenvalues of the corresponding coefficient matrices.
For convenience, we shall often write the eigenvalues of a given matrix $B\in\mathbb{R}^{d\times d}$ 
as $\lambda_{B}$ from now on. 

\begin{lemma}[Matignon \cite{Matignon}] \label{lem:Matignon}
Consider the F-ODEs \eqref{model} with $f\equiv0$. Then it holds that 

\emph{(i)} The solution to \eqref{model} is
 asymptotically stable if and only if all the eigenvalues $\lambda_{A}$ satisfy that 
 \begin{equation} \label{eigen}
\begin{split}
 \lambda_{A}\in \Lambda_{\alpha}^{s}:=\left\{z\in\mathbb{C}\setminus \{0\}: |\arg(z)|> \frac{\alpha\pi}{2}\right\}. 
 \end{split}
\end{equation}

\emph{(ii)} If all eigenvalues $\lambda_{A}\in \Lambda_{\alpha}^{s}$, 
the solution to \eqref{model} is Mittag-Leffler stable, i.e., $\|y(t)\|= O(t^{-\alpha})$ as $t\to \infty$.
\end{lemma}

The result in Part (i) can be found in the appendix of \cite{Lubich86}. 
In order to avoid possible technical complications, 
we shall not consider the critical case in this work, i.e., $|\arg(\lambda_{A})|={\alpha\pi}/{2}$.
Compared with the classical linear ODEs,  we can see two major differences for the F-ODEs:
the system can be still stable when the real parts of the eigenvalues of the coefficient matrix are 
strictly positive (as long as $\lambda_{A}\in\Lambda_{\alpha}^{s}$); 
the long-time polynomial decay rate of the solutions is generally 
slower than the exponential decay rate of the classical ODEs, 
which reflects the non-local nature of fractional derivatives in some sense and 
is also the main motivation for time fractional differential equations to apply in many practical models
describing various slow processes, such as anomalous diffusion \cite{Pod}.

There have been many recent studies about the stability and decay rates of time fractional equations \cite{Cong18, Cong19, Cuesta07, LiCP}.
Cong et al. established the general theory of asymptotic stability of solutions to F-ODEs with constant coefficients 
\cite{Cong18, Cong19} under small perturbations.
The long time behavior of the solutions of the time fractional PDEs can be found in \cite{Zacher19} and the references therein. 
Other relevant stability results for linear or nonlinear time fractional equations, 
with or without time-delay, can be found in the review paper \cite{LiCP}.

For the general F-ODEs \eqref{model}, the nonlinear term $f(t, y)$ can be seen as some perturbation of the corresponding linear system if  $f(t, y)$ is small in some sense. As for classical ODEs, it is natural to expect that if the original unperturbed systems are asymptotically stable and the perturbations $f(y)$ is small enough, then the perturbed systems are also asymptotically stable. We can even guess that 
the trivial solutions of the perturbed systems also have a polynomial decay rate 
under certain conditions, similarly to the original unperturbed systems.  Following such a path, the rigorous perturbation theory of F-ODEs was recently established in  \cite{Cong18,Cong19}, by combining the fractional Lyapunov method and the 
fixed-point technology of Lyapunov-Perron operator.

\begin{lemma} \label{lem:Cong}
(\cite{Cong18,Cong19}) Assume that the spectrum of constant coefficient matrix $A$ satisfies that $\lambda_{A}\in \Lambda_{\alpha}^{s}$. 

\emph{(a)} Assume that the nonlinear perturbation $f(t, y)$ satisfies that 
\begin{equation} \label{eq:fcond}
\begin{split}
f(t,0)=0,\quad~\|f(t,x)-f(t,y)\|\leq L(t)\|x-y\|, \quad \forall\,t\geq 0 \mbox{~~and} ~~x, y\in\mathbb{R}^{d},
\end{split}
\end{equation}
where $L(t): [0, \infty)\to \mathbb{R}_{+}$ is a continuous Lipschitz function and satisfies one of the three conditions:
\begin{equation} \label{eq:Lcond}
\begin{split}
\emph{(i)} &~q_1:=\sup_{t\geq 0} \int_{0}^{t} (t-s)^{\alpha-1}\|E_{\alpha, \alpha} ((t-s)^{\alpha}A)\| L(s)ds<1, 
\quad(robust~ stability);\\
\emph{(ii)} &~ L(t)<q_2 := \sup_{t\geq 0} \frac{1}{2\int_{0}^{t} t^{\alpha-1}\|E_{\alpha, \alpha} (t^{\alpha}A)\|dt}, \quad (uniform~ small~ perturbation);\\
\emph{(iii)} &~ \lim_{t\to \infty}L(t)=0, \quad (decaying~ perturbation).
\end{split}
\end{equation}
Then the trivial solutions to F-ODEs \eqref{model} is asymptotical stable. 

\emph{(b)} Assume that the perturbation $f(y)$ is Lipschitz continuous in a neighborhood of the origin such that 
\begin{equation} \label{eq:fd2}
\begin{split}
f(0)=0\quad \text{and}\quad \lim_{r\to 0} \ell_{f}(r)=0,\quad \text{where}\quad~\ell_{f}(r):=\sup_{x, y\in B(0, r)} \frac{\|f(x)-f(y)\|} {\|x-y\|}.
\end{split}
\end{equation}
Then the trivial solution is Mittag-Leffler stable with optimal decay rate, i.e.,  $\|y(t)\|=O(t^{-\alpha})$ as $t\to\infty$. 
\end{lemma}

Note that the results in \cref{lem:Cong}(b) still hold if we replace the condition $\lim_{r\to 0} \ell_{f}(r)=0$ in \eqref{eq:fd2} with the slightly stronger assumptions $f(y)\in C^{1}$ and $f'(0)=0$.
A concrete example of Mittag-Leffler stability of the solutions for fractional SIRS epidemic model was recently analyzed 
in \cite{MLexam}.

The main tasks of this work are to establish the numerical Mittag-Leffler stability for general F-ODEs in $\mathbb{R}^{d} 
(d\geq 1)$ with or without small perturbations, and to derive the same long-time polynomial decay rate of the numerical solutions 
as the one of the solutions to the continuous equations (cf.\,Lemmas\,\ref{lem:Matignon} and \ref{lem:Cong}), 
as described in detail below:

(1) For homogenous F-ODEs, we consider the numerical Mittag-Leffler stability for F-LMMs and $\mathcal{L}$1 method. 
The main ingredients are the numerical 
stable region characterized by generating functions \cite{Lubich86} and the singularity analysis for generating functions developed in \cite{FS09}. Our analysis helps establish the optimal polynomial long-time behavior of numerical solutions for $\mathcal{L}$1 method and strong $A$-stable F-LMMs, for which the generating polynomial $\delta(z)$ has no poles or zeros in the neighborhood of the unit disk with the exception $z=1$ (see \eqref{eq:Conw}).  This avoids the special requirement by 
the energy method \cite{Xu16,Wang18,WangZou19} for the signs of the coefficients in numerical schemes, 
which would exclude the F-BDF$2$ (F-BDF$k$ is referred to as the F-LMMs or convolution quadratures 
generated by $k$-step BDF in this work).
We shall prove that F-BDF$2$ has an optimal decay rate for homogenous F-ODEs without any stepsize constraint.  
Another advantage of our analysis is its generosity for a unified framework for studying the Mittag-Leffler stability through generating functions. The asymptotic behavior was obtained in \cite{Cermak15} 
for Gr\"{u}nwald-Letnikov method (i.e., F-BDF$1$), but the method there depends heavily 
on the simple and special structure of the coefficients and does not apply to other numerical schemes
such as $\mathcal{L}$1 and F-BDF$k$ for $k\geq 2$.

(2) For non-homogenous F-ODEs with small perturbations, we consider two types of numerical methods. The first is the strong $A$-stable F-LMMs. We first express the numerical solutions as a discrete constant variation formula. 
Then the key step is to derive the asymptotic decay rate of the discrete fractional resolvent operator. We point out that neither the existing singularity analysis for generating functions nor the standard resolvent estimate can help achieve the desired results. Instead, a new estimate based on Prabhakar function is provided to obtain the optimal decay rate for the discrete resolvent operator, which is exactly consistent with the continuous case. Using this result, we are able to establish the numerical Mittag-Leffler stability for strong $A$-stable F-LMMs for F-ODEs with small perturbations.

The second type of numerical schemes are those 
based on a $\alpha$-difference approximation to Caputo derivatives, 
which can be seen as a fractional extension of the backward Euler formula.  
Numerical solutions to these schemes can also be expressed by the discrete 
constant variation formula involving a discrete fractional resolvent operator.  A significant advantage of this approach is that the discrete fractional resolvent operator can be directly connected with the continuous one through the Poisson transformation \cite{Lizama17, Ponce20}. Therefore, we can use various properties of Poisson transformation and continuous operator to give an optimal estimate of the decay rate for the discrete fractional resolvent operator in a simple and straightforward manner. Then combined with our earlier results about the asymptotic behavior of the Volterra difference equation \cite{Wang18}, 
we derive the Mittag-Leffler stability of numerical solutions under the natural smallness of the perturbations.

We now recall some existing numerical methods and their stabilities 
for time fractional differential equations. Most schemes may be classified in two major groups, namely  
F-LMMs developed by Lubich in 1980s \cite{Lubich85, Lubich86, Lubich88} 
and interpolation based methods (such as $\mathcal{L}$1-type methods).
F-LMMs have become very popular because they inherit the good stability of classic LMMs 
and can be implemented easily. 
Note that F-LMMs can be interpreted as convolution quadratures  in a more general way, 
which is an efficient numerical approach for approximating general convolution operators  \cite{Lubich88}. 
The $\mathcal{L}$1 schemes are among the most popular numerical approximations for Caputo derivatives, and are 
easy to implement with acceptable precision and good numerical stability. 
This method was systematically studied in \cite{Jin16, Yan18} for sub-diffusion equations.
Since the $\mathcal{L}$1 scheme approximates 
the classical first derivative by backward Euler scheme on each subinterval,
it provides a good basis for numerical approximations for Caputo derivatives on non-uniform grids \cite{Kopteva19,Liao19}.
We refer to the recent survey \cite{Martin21} for various properties and applications of the $\mathcal{L}$1 method. 
We point out here that the method in this paper cannot be applied to non-uniform grids at present. The main reason is that for non-uniform grids, time step size is no longer a constant and the weight coefficients depend on step size, so it is difficult to define the related generating function.

One of the main difficulties in solving time fractional differential equations numerically is the limited regularity of the true solutions near the initial time at $t=t_{0}$.
The limited regularity often causes some reduction of convergence rates of numerical schemes. 
Special correction techniques or specific non-uniform meshes can be developed near the initial time 
to restore optimal convergence rates of numerical schemes 
for time fractional evolutional equations 
\cite{Jin17, Kopteva19, Liao19,Lubich85,Martin21,Yan18}.

In this work, we are mainly concerned with an important mathematical issue whether 
the numerical solutions can inherit the long-term qualitative characteristics of 
the solutions to the continuous problems. 
As Lemmas \ref{lem:Matignon} and \ref{lem:Cong} indicated, the Mittag-Leffler stability with long-time polynomial decay rate of solutions is a key characteristic of the F-ODEs. As far as we know,  
no much study exists on the qualitative behaviors of numerical solutions to time fractional differential equations.  
Cuesta et.al. established the asymptotic behavior for both continuous and discrete solutions
to an abstract time fractional ingegro-differential equations of order $\alpha\in (1, 2)$ 
\cite{Cuesta03,Cuesta07}.
We analyzed recently in \cite{Wang18} the contractivity, 
dissipativity and long-time polynomial decay rate of solutions to 
Gr\"{u}nwald-Letnikov formula and $\mathcal{L}$1 method for a class of real nonlinear F-ODEs.   
It is extended to more complex stiff time fractional functional differential equations in \cite{WangZou19}. 
However, those analyses can not apply to study the long-time behavior of numerical solutions to 
F-ODEs \eqref{model}, even for the simple case $f\equiv0$, 
because the eigenvalues $\lambda_{A}$ can be complex or have positive real parts. 
Another disadvantage of the energy methods used in \cite{Wang18} often results in a suboptimal decay rate, 
i.e., $\|y_{n}\|=O( t_{n}^{-\alpha/2})$ rather than $\|y_{n}\|=O(t_{n}^{-\alpha})$; 
see detailed explanations in \cite{Wang18}.  
For linear evolutionary Volterra integro-differential equations in Hilbert space, the uniform behavior of numerical methods were derived in \cite{Xu19} and the references therein.

The rest of the paper is organized as follows.   
In \cref{sec:F-LMMs},  we recall a singularity analysis of generating functions and some basic concepts and properties of F-LMMs and $\mathcal{L}$1 schemes, and then slightly generalize the stability results to the general vector F-ODEs.
In \cref{sec:homogenous}, we establish the numerical Mittag-Leffler stability for both F-LMMs and $\mathcal{L}$1 schemes 
for homogeneous linear vector F-ODEs.  In \cref{sec:Perturbations}, we consider the numerical Mittag-Leffler stability for F-ODEs with small perturbations for F-LMMs and a numerical scheme based on $\alpha$-difference approximation to Caputo fractional derivative, respectively. 
We present in \cref{sec:examples} several typical numerical examples to illustrate and verify our theoretical results. 

In the subsequent analysis, we will often use $C$ to represent a generic positive constant, 
which may take different values at different occasions but is always independent of $t$ and $t_n$.
We will also write the convolutional identity as $\delta_{d}:=(1, 0, 0,...)$, 
and $\delta_{i, j}$ as the Kronecker function, i.e., 
$\delta_{i, j}=1$ if $i=j$ and $\delta_{i, j}=0$ if $i\neq j$, and 
$\delta_{n,0}$ is the $n$th entry of $\delta_{d}.$
Furthermore, we shall often write $F_{v}(z)\sim f(z)$ as $z\to z_{0}$ to stand for their equivalence 
in the sense ${F_{v}(z)}/{f(z)}\to 1$ as ${z\to z_{0}}$.

\section{Numerical methods for F-ODEs}
\label{sec:F-LMMs}
We start with some basic concepts and notation.
We write a discrete sequence as $v=(v_{0}, v_{1},...)$ for $v_{n}\in\mathbb{C}^{d}$. If $u, v$ are two scalar sequences with $u_{n}, v_{n}\in\mathbb{C}^{1}$, we define the discrete convolution $u*v=w$, with $w_{n}=\sum_{j=0}^{n}u_{n-j}v_{j}$. When $u*v=\delta_{d}$, 
we say the sequence $u$ is invertible, and write the inverse $v=u^{(-1)}$. Note that a sequence $u$ is invertible if and only if $u_{0}\neq 0$. We often write $[\cdot]_{n}$ to be the $n$-th entry of a sequence, i.e., $[u]_{n}=u_{n}$.  The generating function of a sequence $v=(v_0, v_1, \ldots)$ is defined by 
$
F_v(z)=\sum_{n=0}^{\infty} v_n z^n, z\in\mathbb{C}. 
$
It is easy to verify that $F_{u*v}(z)=F_u(z)F_v(z)$. Hence, the generating functions of $\mu$ and $\omega$ are related by
$F_{\mu}(z)={1}/{F_{\omega}(z)}$ if $\omega=\mu^{(-1)}$.

Consider the F-ODE $\mathcal{D}_t^{\alpha}y(t) =g(t, y(t))$, subject to $y(0)=y_{0}$.
The implicit scheme approximating $y(t_{n})$ by $y_{n}$ ($n\geq 1$) at the uniform grids $t_{n}=nh$ with step size $h>0$ has the following form:
\begin{gather}\label{eq:Volterradiff}
\mathcal{D}_h^{\alpha}(y_{n}) :=\frac{1}{h^{\alpha}}
\sum_{j=0}^n \mu_{j}(y_{n-j}-y_0)=g(t_{n}, y_n):=g_n,~~n\ge 1.
\end{gather}
The operator $\mathcal{D}_h^{\alpha}(y_{n})$ is the numerical approximation of the Caputo derivative $\mathcal{D}_t^{\alpha}y(t)$ at $t=t_{n}$ and the form in \eqref{eq:Volterradiff} can be viewed as a fractional backward differential formula. The weight coefficients $\{\mu_{j}\}_{j=0}^{\infty}$ can be determined in several ways. If we like to include the case $n=0$, the equation \eqref{eq:Volterradiff} is written as
\begin{gather}\label{eq:Voldiff}
\frac{1}{h^{\alpha}} \sum_{j=0}^n \mu_{j}(y_{n-j}-y_0)= g_n-g_0\delta_{n,0},~~n\ge 0.
\end{gather}

It is known that the F-ODEs are equivalent to the Volterra integral equations of second class, and numerical schemes can be constructed from the integral form. In fact, by using the discrete convolution inverse, we can write the numerical scheme \eqref{eq:Voldiff} 
as the equivalent integral form:
\begin{gather}\label{eq:Volinter}
\begin{split}
y_n-y_0&=h^{\alpha}[\omega*(g- g_0 \delta_{d})]_n=h^{\alpha}[\omega*g-g_0 \omega]_{n}=h^{\alpha}\sum_{j=1}^{n} \omega_{n-j} g_{j}, ~~n\ge 1.\\
\end{split}
\end{gather}
We shall always use $\omega=(\omega_0, \omega_1,...)$ and $\mu=(\mu_0, \mu_1,...)$ 
for the coefficients in the integral form and differential form, respectively, and both are related by $\omega=\mu^{(-1)}$.

Compared with the continuous F-ODEs, the term in the right-hand side of \eqref{eq:Volinter} corresponds to integral approximation of Riemann-Liouville fractional integral. The weight coefficients $\{\omega_{j}\}_{j=0}^{\infty}$ will directly link to the stand kernel {$k_{\alpha}(t):={t^{\alpha-1}}/{\Gamma(\alpha)},\, t>0$.  
A class of $\mathcal{CM}$-preserving} numerical schemes were developed in \cite{LiWang19}, 
by means of {the complete monotonicity} of $k_{\alpha}(t)$.
In this paper, we focus on another key property of the solutions of F-ODEs, the Mittag-Leffler stability, which highlights 
the long-term optimal polynomial decay rate of the solutions. 
The generating function will be a main tool in our analysis since
the asymptotic properties of a sequence can be characterized in terms of its generating function.

\begin{lemma}\emph{(\cite[Corollary VI.I]{FS09})} \label{lem:generating}
 Assume $F_v(z)$ is analytic on $\Delta(R,\theta) :=\{z: |z|<R, z\neq 1, |\mathrm{arg}(z-1)|>\theta\}$ for some $R>1$ and $\theta\in (0, \frac{\pi}{2})$. If
 $F_v(z)\sim (1-z)^{-\beta}$ as $z\to 1, z\in \Delta(R,\theta)$ for $\beta\neq \{0, -1, -2, -3, \cdots\}$, then $v_n \sim \frac{1}{\Gamma(\beta)}n^{\beta-1} $ as $n\to \infty$.
\end{lemma}

The singularity analysis of generating function was  
developed in \cite{FO90}, and further extended to various typical functions in the monograph \cite{FS09}.
The above fundamental lemma describes the correspondence between the asymptotic expansion of a function near its dominant singularities and 
the asymptotic expansion of the function's coefficients.  The advantages of this approach is that it only necessitates local asymptotic properties of the generating functions, and hence can be applied to functions whose singular expansions contain fractional powers and is well suited for the generating functions of F-LMMs and $\mathcal{L}$1 schemes.

\subsection{F-LMMs}
In a series of pioneering work \cite{Lubich85,Lubich86, Lubich88}, Lubich proposed the fractional linear multistep methods (F-LMMs) for weakly singular Abel-Volterra integral equations, which can be directly applied to Caputo F-ODEs. 
Let us first recall the basic concepts and results of this approach, and 
consider the numerical approximating of Riemann-Liouville fractional integral 
\begin{equation} \label{eq:RL}
\mathcal{I}_{t}^{\alpha} y(t)= \int_{0}^{t}k_{\alpha}(s) y(t-s)ds
= \int_{0}^{t} \frac{1}{2\pi i} \int_{\mathcal{C}} e^{s \lambda}K_{\alpha}(\lambda) y(t-s)d\lambda ds
=\frac{1}{2\pi i} \int_{\mathcal{C}}  \int_{0}^{t}e^{s \lambda}y(t-s)  K_{\alpha}(\lambda) ds d\lambda,
\end{equation}
where $K_{\alpha}(\lambda)=\lambda^{-\alpha}$ 
is the Laplace transform of the standard kernel $k_{\alpha}(t)$, 
and the contour $\mathcal{C}$ is properly chosen for the well-definedness of the integral.   
We can easily see that $u(t):= \int_{0}^{t}e^{s \lambda}y(t-s)ds$ solves the ODE
$u'(t)=\lambda u(t)+y(t)$ for $t>0$, with the initial value  $u(0)=0$.
We apply the classical $k$-step LMMs with generating polynomials $\rho(z)=\sum_{j=0}^{k}\alpha_{j}z^j$ and  $\sigma(z)=\sum_{j=0}^{k}\beta_{j}z^j$ to the ODE to obtain the numerical solution $u_{n}$ approximating $u(t_n)$.
By putting the initial values $u_{-k}=...=u_{-1}=0$ and introducing the generating functions $F_{u}(z)=\sum_{j=0}^{\infty}u_{j}z^{j}$ and $F_{y}(z)=\sum_{j=0}^{\infty}y_{j}z^{j}$, one can verify that $F_u(z)$ satisfies the equation 
\begin{equation} \label{eq:kGF}
\begin{split}
\left(\alpha_{0}z^k+\cdots+\alpha_{k-1}z+\alpha_{k} \right)F_u(z)=\left( \beta_{0}z^k+\cdots+\beta_{k-1}z+\beta_{k} \right)(h\lambda \cdot F_u(z)+hF_y(z)).
\end{split} 
\end{equation}
Solving this equation for $F_u(z)$ yields that  
\begin{gather}\label{eq:delta}
\begin{split}
u_n= \left[ \left(\frac{\delta(z)}{h}-\lambda \right)^{-1} F_{y}(z)\right]_{n}, \text{ where }
\delta(z)=\frac{z^k \rho(z^{-1})}{z^k \sigma(z^{-1})}=\frac{\alpha_{0}z^k+\cdots+\alpha_{k-1}z+\alpha_{k}} {\beta_{0}z^k+\cdots+\beta_{k-1}z+\beta_{k}}.\\
\end{split}
\end{gather}
Substituting this formula into \eqref{eq:RL} and applying the Cauchy integral formula, we arrive at 
\begin{equation} \label{eq:RLintegral}
\begin{split}
\mathcal{I}_{t_{n}}^{\alpha} y(t_{n})&=\frac{1}{2\pi i} \int_{\mathcal{C}}  \left( \int_{0}^{t_n}e^{s\lambda}y(t_n-s)ds\right) K_{\alpha}(\lambda) d\lambda 
\approx \frac{1}{2\pi i} \int_{\mathcal{C}}  \left[  \left(\frac{\delta(z)}{h}-\lambda \right)^{-1} F_{y}(z) \right]_{n} K_{\alpha}(\lambda) d\lambda\\
&=\left[  K_{\alpha}\left(\frac{\delta(z)}{h} \right) F_{y}(z) \right]_{n}=\left[  \left(\frac{\delta(z)}{h} \right)^{-\alpha} F_{y}(z) \right]_{n}
:=h^{\alpha}\sum_{j=0}^{n}\omega_{n-j} y_{j} 
\end{split} 
\end{equation}
for $n\ge 1$, where we have written 
\begin{equation} \label{eq:Conw}
\begin{split}
F_{\omega}(z):=(\delta(z))^{-\alpha}=\sum_{j=0}^{\infty}\omega_j z^{j}.
\end{split} 
\end{equation}

Now we can readily get the F-LMMs for F-ODEs \eqref{model} 
by using the approximation \eqref{eq:RLintegral} in 
the corresponding Volterra integral equations:
\begin{gather}\label{eq:FLMMs0}
y_n=y_0+h^{\alpha}\sum_{j=0}^{n} \omega_{n-j} \left(Ay_{j}+f(t_j, y_{j})\right), ~~n\ge 1.
\end{gather}
However, we can see the two numerical schemes \eqref{eq:Volinter} and  \eqref{eq:FLMMs0} are different. They differ by a term related to the initial value $h^{\alpha} \omega_{n} g_{0}$, with $g=Ay+f(t, y)$.
In order for the two methods to be consistent, we can make some appropriate modifications to \eqref{eq:FLMMs0}. 
An easy way to do this is to drop the term for $j=0$ in \eqref{eq:FLMMs0}, which leads to the modified approximation 
$
\mathcal{I}_{t_{n}}^{\alpha} g(t_{n})\approx h^{\alpha}\sum_{j=1}^{n}\omega_{n-j} g_{j}.  
$
In fact, this correction method was used in \cite[pp.4, Eq. (1.15)]{Lubich96} for F-BDF$1$ to obtain a positive definite discrete quadrature formula.
Through this correction, we obtain the F-LMMs for F-ODEs: 
\begin{gather}\label{eq:FLMMs}
y_n=y_0+h^{\alpha}\sum_{j=1}^{n} \omega_{n-j} \left(Ay_{j}+f(t_j, y_{j})\right), ~~n\ge 1\,, 
\end{gather}
where $\omega_j$ are still given by \eqref{eq:Conw}. The numerical method \eqref{eq:FLMMs} is now also 
consistent at $n=0$ by noting the sum is zero when the upper index is smaller than the lower index,
and also fully consistent with the scheme derived in \eqref{eq:Volinter} by the convolution inverse.


\subsection{Stability region of F-LMMs} 
For classical ODEs (i.e., $\alpha=1$), we know that the linear test equation $y'(t)=\lambda y(t)$
is asymptotic stable if and only if $\mathrm{Re}(\lambda)<0$. So a reasonable numerical scheme is 
often required to preserve this stability \cite{Hairer}.
On the other hand, the linear test model for F-ODEs (with $\alpha\in(0,1)$) is 
$\mathcal{D}_{t}^{\alpha}y(t)=\lambda y$, which is 
asymptotic stable if $\lambda\in \Lambda_{\alpha}^{s}$, as indicated in \cref{lem:Matignon}. 
$\Lambda_{\alpha}^{s}$ coincides with the left complex semi-plane when $\alpha=1$.
Applying F-LMMs to the fractional linear test equation gives 
\begin{gather}\label{eq:linearstability}
y_n=y_0+\lambda h^{\alpha} [\omega*(y-y_0\delta_d)]_n,~~n\ge 0,
\end{gather}
where the weights $\omega_{n}$ are given in \eqref{eq:Conw}.
The numerical stability region for \eqref{eq:linearstability} is defined by  
\begin{gather}\label{eq:numstab}
\mathcal{S}_h^{\alpha}:=\{\zeta=\lambda h^{\alpha}\in \mathbb{C}\setminus \{0\}:  y_n\to 0~ \hbox{as}~ n\to \infty \},
\end{gather}
and a numerical method is said to be $A(\beta)$-stable (with $\beta\in(0,\pi)$), if the region $\mathcal{S}_h^{\alpha}$ contains the infinite wedge 
$
A(\beta)=\{z \in\mathbb{C}\setminus \{0\}; |\arg(-z)|< \beta \}.
$
Similarly to ODEs,
if the stability region $\mathcal{S}_{h}^{\alpha}$ contains the entire sector $\Lambda_{\alpha}^{s}$,
i.e., $\mathcal{S}_{h}^{\alpha} \supset \Lambda_{\alpha}^{s}$,
then the method is said to be $A(\frac{\alpha\pi}{2})$-stable, or simply $A$-stable. 

Let us now recall the definition of the strong $A(\beta)$-stability of a LMM 
defined by a generating polynomial
$F_{\overline{\omega}}(z)=F_{\overline{\omega}(\rho, \sigma)}(z)=\delta(z)^{-1}$
for the classical ODE, with order $p\geq 1$ \cite{Lubich85, Lubich86,  Lubich88}: 
 \begin{equation} \label{eq:assum}
\begin{split}
 \qquad \begin{array}{l}\delta(z) \text { is analytic, with no zeros in a neighborhood of the  
unit disk $|z| \leq 1$ except } z=1; \\ 
|\arg \delta(z)| \leq \pi-\beta \text { for }|z|<1; 
~~\frac{1}{h} \delta\left(\mathrm{e}^{-h}\right)=1+O\left(h^{p}\right), \text { with } p \geq 1\,,
\end{array}
\end{split}
\end{equation}
where $\delta(z)$ is defined in \eqref{eq:delta}. 
We note that the conditions \eqref{eq:assum} exclude the simple trapezoidal rule (i.e., $\delta(z)=2(1-z)/(1+z)$), 
which is $A$-stable but not strongly $A$-stable.
The following lemma presents the fundamental relationship 
between the stability regions of the classical LMMs and the F-LMMs.

\begin{lemma}
\label{lem:Lubich}
(\cite{Lubich85}) Consider a classical LMM defined by a generating polynomial 
$F_{\overline{\omega}}(z)=\delta(z)^{-1}$ satisfies the stability conditions \eqref{eq:assum}.
Let $\mathcal{S}_{h}$ and $\mathcal{S}_{h}^{\alpha}$ be 
the stability regions of the standard LMM  and its corresponding 
F-LMM defined by $F_\omega(z)=\left(F_{\overline{\omega}}(z)\right)^{\alpha}=\delta(z)^{-\alpha}$ respectively.
Then it holds that 

\emph{(i)}  $\mathcal{S}_{h}^{\alpha}=\mathbb{C}\setminus \left\{1/F_{\omega}(z): |z|\le 1\right\}$;

\emph{(ii)} $\left(\mathbb{C}\setminus \mathcal{S}_{h}^{\alpha}\right)=\left(\mathbb{C}\setminus \mathcal{S}_{h}\right)^{\alpha}$;

\emph{(iii)} LMM is $A$-stable if and only if the F-LMM is $A$-stable;

\emph{(iv)} with $\pi-\varphi=\alpha(\pi-\psi)$, LMM is $A(\varphi)$-stable if and only if the F-LMM is $A(\psi)$-stable.

\end{lemma}

\subsection{Examples of F-LMMs} 
Let us recall some typical examples of F-LMMs, which are direct extensions of LMMs \cite{Hairer} and 
inherit their good numerical stability due to \cref{lem:Lubich}.

\subsubsection{F-BDF$k$} 
The generating functions of F-BDF$k$ are given by 
 \begin{equation} \label{eq:FBDFk}
\begin{split}
F_{\mu}(z):=\frac{1}{F_{\omega}(z)}=\left(\sum_{\ell=1}^{k}\frac{1}{\ell}(1-z)^{\ell}\right)^{\alpha}
=\sum_{j=0}^{\infty}\mu_{j}z^{j},~~~k=1,2,...,6.
\end{split}
\end{equation}

For $k=1$, we have $F_{\mu}(z):=(1-z)^{\alpha}$, and the scheme is just the Gr\"{u}nwald-Letnikov formula.
The weights $\mu_{j}$ can be recursively evaluated as
$\mu_{0}=1$ and $\mu_{j}=\left(1-\frac{\alpha+1}{j}\right)\mu_{j-1}$ for $j\geq 1$.

For $k=2$, we have $F_{\mu}(z):=\left(\frac{3}{2}-2z+\frac{1}{2}z^{2} \right)^{\alpha}=\sum_{n=0}^{\infty}\mu_{j}z^{j}$, 
and the weights $\{\mu_{j}\}_{j=0}^{\infty}$ satisfy that
 \begin{equation*}
\begin{split}
&\mu_{0}=\left(\frac{3}{2}\right)^{\alpha}>0, ~\mu_{1}=-\left(\frac{3}{2}\right)^{\alpha}\frac{4\alpha}{3}<0, ~\mu_{2}=\left(\frac{3}{2}\right)^{\alpha}\frac{\alpha(8\alpha-5)}{9},\\
&\mu_{3}=\left(\frac{3}{2}\right)^{\alpha}\frac{4\alpha(\alpha-1)(7-8\alpha)}{81},~ \mu_{j}<0~\text{for}~ j\geq 4,~~ \text{and~that}~\sum_{j=0}^{\infty}\mu_{j}=0.\\
\end{split}
\end{equation*}

Both F-BDF$1$ and F-BDF$2$ are $A$-stable. One important characteristic of F-BDF$2$ is that the coefficients $\mu_{j}$ for $j\geq 1$ are not all negative, e.g., $\mu_2$ and $\mu_3$ could be positive. This is very different 
from F-BDF$1$, whose coefficients are all negative. {From \cite{LiWang19}, we know F-BDF$1$ is $\mathcal{CM}$-preserving but F-BDF$2$ is not. }
In fact, it is a common feature that the coefficients of higher order methods for time fractional derivatives do not keep 
same sign. This often causes some difficulties in convergence analysis when energy-type methods are used; 
see \cite{Xu16, Wang18,WangZou19}.
As we see from our subsequent results, it is unnecessary to impose any step size requirements on 
the $A$-stable F-BDF$2$ to preserve the long-time polynomial decay rate of the solutions to 
linear time fractional evolutional equations.

\subsubsection{F-Adams schemes} 
$k$-step fractional Adams methods are generated by the generating function
 \begin{equation} \label{eq:FAdams}
\begin{split}
F_{\omega}(z)=(1-z)^{-\alpha} \left(\gamma_{0}+\gamma_{1}(1-z)+...+\gamma_{k}(1-z)^{k} \right),
\end{split}
\end{equation}
where the parameters $\gamma_{j}$ are the coefficients in the truncated expansion of the $\alpha$-power 
of function $G(t)=\frac{-t}{\ln(1-t)}$:
\[
(G(1-z))^\alpha=\left(\frac{1-z}{-\ln(z)} \right)^\alpha= \sum_{j=0}^{\infty} \gamma_{j} (1-z)^{j},
\]
which is used to generate the Adams-Moulton methods for ODEs.
The consistency and convergence of order $k$ of the F-Adams scheme were proved in \cite{Lubich86}.
It is easy to see that $\gamma_{0}=1, \gamma_{1}=-\frac{\alpha}{2}$. 
We note that the $1$-step method is just the Gr\"{u}nwald-Letnikov scheme, while 
the $2$-step second order fractional Adams method has its generating function 
 \begin{equation}\label{FAdams}
F_{\omega}(z)=(1-z)^{-\alpha} \left(1-\frac{\alpha}{2}(1-z) \right).
\end{equation}
This method is also $A$-stable. Fast algorithms for computing the weights of the fractional power $(1\pm z)^{\alpha}$ based on Miller formula \cite{Garrappa15} can be used to evaluate the expansion coefficients of most $A$-stable F-LMMs effectively.

\subsection{$\mathcal{L}$1 method and its numerical stability region}
The $\mathcal{L}$1 scheme can be seen another fractional generalization of the backward Euler scheme for ODEs \cite{Martin21}. 
The $\mathcal{L}$1 scheme approximating the Caputo fractional derivative can be written in the discrete convolution form
\[
\mathcal{D}_h^{\alpha}(y_n):=\frac{1}{h^{\alpha}} \left(\sum_{j=0}^{n-1} \mu_{j}y_{n-j}-\sigma_n y_0 \right)
= \frac{1}{h^{\alpha}} \sum_{j=0}^{n} \mu_{j}(y_{n-j}-y_0),
\] 
where
$\mu_{0}=\frac{1}{\Gamma(2-\alpha)}, \sigma_{n}= \frac{1}{\Gamma(2-\alpha)} \left( n^{1-\alpha}-(n-1)^{1-\alpha}\right)$ and
$\mu_{j}=\frac{1}{\Gamma(2-\alpha)} \left( (j+1)^{1-\alpha}-2j^{1-\alpha}+(j-1)^{1-\alpha}\right)$  for $j\geq 1.$
The generating function of the $\mathcal{L}$1 scheme is given by 
\begin{gather}\label{eq:generatingfunL1}
F_\mu(z)=\sum_{n=0}^{\infty} \mu_n z^n= \frac{1}{\Gamma(2-\alpha)} \left(\frac{1}{z}-2+z\right) \mathrm{Li}_{\alpha-1}(z), 
\end{gather}
where $\mathrm{Li}_{\gamma}(z)=\sum_{n=1}^{\infty}{z^n}/{n^\gamma}$ is the polylogarithm function.
The rigorous stability analysis of the $\mathcal{L}$1 scheme is more difficult than that of F-LMMs 
due to the involvement of the polylogarithm function. 
It was proved in \cite{Jin16} that $\mathcal{L}$1 scheme is at least $A(\pi/4)$-stable 
and the result can be improved to be the $A$-stable 
by making use of a very elaborate expansion formula of the polylogarithm function.

\subsection{Stability region for vector-valued F-ODEs}
We now extend the stability results in \cref{lem:Lubich} for F-LMMs from the scalar test equation 
to the general vector-valued F-ODEs,
which are used to prove the discrete version of the stability result in \cref{lem:Matignon}.  
The key in the proof of \cref{lem:Lubich} is the application of the discrete Paley-Wiener theorem and 
a technique to deal with the singularity of the generating function $F_{\omega} (z)$ at $z=1$ \cite{Lubich85, Lubich86}. 
We follow this idea and apply the following vector-valued version of the discrete Paley-Wiener theorem \cite{Lubich85}.

\begin{lemma} \label{lem:Lubich2}
Consider the discrete Volterra integral equation 
$y_{n}= p_{n}+\sum_{j=0}^{n}Q_{n-j}y_{j}$ for $n\geq0,$
 where the matrix sequence $\{Q_{n}\}_{n=0}^{\infty}$ belongs to $\ell^{1}$ (i.e., each entry in the sequence is in $\ell^{1}$).
Then it holds that $\|y_{n}\|\to 0$  (resp. bounded) whenever $\|p_{n}\|\to 0$  (resp. bounded) as  $n\to\infty$
if and only if the Paley-Wiener condition is satisfied, i.e.,
\begin{equation}\label{eq:PWcondition}
\det \Big( I- \sum\limits_{j=0}^{\infty}Q_{j}z^{j}\Big)\ne 0 ~~for ~~ |z|\le 1.
\end{equation}
\end{lemma}

\begin{theorem}\label{thm:stability}
Assume that the F-LMM satisfies the conditions in \eqref{eq:assum}. 
Then for the homogenous vector-valued F-ODEs in \eqref{model} with $f\equiv0$ and any $h>0$, 
the numerical stability region is given by
\begin{equation} \label{eq:stabregion}
\begin{split}
\mathbb{S}_{h}^{\alpha}&=\det \left( I- h^{\alpha} F_{\omega} (z)A \right)\ne 0 ~for~ |z|\le 1\\
&\Leftrightarrow \frac{1}{ h^{\alpha} F_{\omega} (z)} \text{~is not an eigenvalue of matrix~} A \text{~for~} |z|\le 1\\
&\Leftrightarrow \mathbb{C}\backslash \left\{\frac{1}{ h^{\alpha} F_{\omega} (z)} \text{~is an eigenvalue of matrix~} A \text{~for~}  |z|\le 1\right\}.
\end{split}
\end{equation}
\end{theorem}

\begin{proof}
A direct application of the F-LMM to the homogenous F-ODEs yields that 
\begin{gather}\label{eq:FLMMH2}
\begin{split}
y_n&=y_0+h^{\alpha}\sum_{j=1}^{n} \omega_{n-j} \left(Ay_{j}\right)
=(I-h^{\alpha}\omega_{n} A)y_{0}+h^{\alpha}\sum_{j=0}^{n}  \left( \omega_{n-j}A\right) y_{j}
\end{split}
\end{gather}
where $\omega_j$ are the coefficients in \eqref{eq:Conw}.
In order to apply the discrete Paley-Wiener theorem, the matrix 
sequence $\{h^{\alpha}\omega_{n}A\}_{n=0}^{\infty}$ is required to be in $\ell^{1}$.
By noting that $A$ is a constant matrix, this condition is equivalent to that the scalar sequence $\{\omega_{n}\}_{n=0}^{\infty}$ is in $\ell^{1}$.

First, the strong stability condition for F-LMMs yields 
$
F_\omega(z)=(1-z)^{-\alpha} u(z),
$
where $u(z)$ is holomorphic in a neighborhood of the unit disc $|z|\leq 1$ \cite[p.467]{Lubich85}.
The main difficulty is now from the fact that the sequence of coefficients in the expansion of $(1-z)^{-\alpha}$ is not in $\ell^{1}$. 
This can be seen from the following expansion for an arbitrary complex number 
$\alpha\in\mathbb{C}\backslash \mathbb{Z}_{\leq 0}$ 
\cite[Theorem VI.I]{FS09}: 
\[
\left[(1-z)^{-\alpha}\right]_{n}\sim \frac{n^{\alpha-1}}{\Gamma(\alpha)} \left(1+\frac{\alpha(\alpha-1)}{2n} +O (n^{-2})\right) \quad \mbox{as } n\to \infty.
\]
Following \cite{Lubich85}, we can add the factor $(1-z)^{\alpha}$ to the generation function 
$F_\omega(z)$ to overcome this difficulty.  
It follows from \eqref{eq:FLMMH2} that $F_{y}(z)=g(z)+h^{\alpha}F_{\omega}(z)AF_{y}(z)$, where $g(z):= \sum_{n=0}^{\infty} (I-h^{\alpha}\omega_{n} A)y_{0}z^{n}$. So if  $\det \left( I- h^{\alpha} F_{\omega} (z)A \right)\ne 0$  for $|z|\le 1$, 
then we have
\begin{gather}\label{eq:Fyfun}
\begin{split}
F_{y}(z)=\left( I-h^{\alpha}F_{\omega}(z)A \right)^{-1} g(z) =\left( (1-z)^{\alpha}I-h^{\alpha}u(z)A \right)^{-1}\cdot  (1-z)^{\alpha}g(z).
\end{split}
\end{gather}
Now we can see that the coefficients of $(1-z)^{\alpha}$ and $u(z)$ are both in $\ell^1$. Hence, Wiener's inversion theorem shows that the coefficients of 
$\left( (1-z)^{\alpha}I-h^{\alpha}u(z)A \right)^{-1}$ is also in $\ell^1$.
On the other hand, it is easy to see that 
\[
(1-z)^{\alpha}g(z)=(1-z)^{\alpha}\left( \frac{1}{1-z}I - h^{\alpha} F_{\omega}(z)A \right)y_{0}= \left( (1-z)^{\alpha-1}I - h^{\alpha} u(z)A \right)y_{0}.
\]
This implies $\|p_{n}\|\to 0$ in the expansion $(1-z)^{\alpha}g(z)=\sum_{n=0}^{\infty}p_{n}z^{n}$.
Now the desired results follow readily from \cref{lem:Lubich2}.
\end{proof}

Theorem\,\ref{thm:stability} indicates that if $\lambda_{A}\in\Lambda^{s}_{\alpha}$ and the F-LMMs are strongly $A$-stable, then $\det \left( I- h^{\alpha} F_{\omega} (z)A \right)\ne 0$ for all $|z|\le 1$ and $h>0$. This means the F-LMMs is unconditionally stable for the vector-valued F-ODEs.  
We can easily see \cref{thm:stability} reduces to \cref{lem:Lubich}(i) for $d=1$.

\section{Mittag-Leffler stability for homogenous F-ODEs} 
\label{sec:homogenous}

In this section, we study the polynomial decay rate of numerical solutions for homogenous linear F-ODEs in 
$\mathbb{R}^d$, 
which can be seen as the discrete version of \cref{lem:Matignon}. 

\begin{theorem}\label{thm:homoFODEs}
Consider the homogenous linear F-ODEs \eqref{model} (i.e., $f\equiv 0$) and assume that  
all the eigenvalues of $A$ satisfy that $\lambda_{A}\in\Lambda_{\alpha}^{s}$.
Then the numerical solutions 
obtained from the strong $A$-stable F-LMMs or $\mathcal{L}$1 scheme 
are Mittag-Leffler stable, i.e., $\|y_{n}\|=O(t_{n}^{-\alpha})$ as $n\to \infty$. 
\end{theorem}

\begin{proof} 
The main idea of the proof is to exploit the special structure of the generating function of F-LMMs and $\mathcal{L}$1 scheme 
so that we can apply \cref{lem:generating}.
We divide the proof into three different cases: 
F-LMM for the scalar test equation, $\mathcal{L}$1 scheme for the scalar test equation, 
and F-LMM and $\mathcal{L}$1 scheme for the general vector-valued system.

{\bf Case I}: F-LMM for the scalar test equation. 
Applying the F-LMM to the scalar test equation, we get
\begin{gather*}
y_n=y_0+\lambda h^{\alpha} [\omega*(y-y_0\delta_d)]_n,~~n\ge 0,
\end{gather*}
where the coefficients $\omega=(\omega_0, \omega_1,\cdots)$ is given by the generating function in \eqref{eq:Conw}.
Multiplying both sides of the equation with $z^{n}$ and summing over $n\ge 0$, we obtain that 
\begin{equation*} 
\begin{split}
F_y(z)=y_0 (1-z)^{-1}+\lambda h^{\alpha}  \left( F_\omega (z) F_{y} (z)-y_0 F_\omega (z) \right).
\end{split}
\end{equation*}
This indicates the generating function for the numerical solution sequence $\{y_n\}$:
\begin{gather}\label{eq:genefun}
F_y(z)=y_0\frac{(1-z)^{-1}-\zeta F_{\omega} (z)}{1-\zeta F_{\omega}(z)}
=y_0 \left(1+\frac{z}{(1-\zeta F_\omega(z))(1-z)}\right),
\end{gather}
where $\zeta=\lambda h^{\alpha}\in\Lambda_{\alpha}^{s}$.
In order for the function $F_{y}(z)$ to be analytic in the region 
\begin{gather*}
\Delta(R,\theta)=\{z\in\mathbb{C}: |z|\le R, z\neq 1, ~|\arg(z-1)|>\theta \}
\end{gather*}
for some $R>1$ and $\theta\in (0, \frac{\pi}{2})$, it is sufficient to require that  
\begin{gather}\label{eq:stabcondition}
1-\zeta F_\omega(z)\neq 0,\quad z\in \Delta(R,\theta).
\end{gather}

By the $A$-stability of F-LMMs, we know $\zeta\in\Lambda_{\alpha}^{s}\subseteq \mathcal{S}_h^{\alpha}$, 
and further obtain $\mathcal{S}_{h}^{\alpha}=\mathbb{C}\setminus \left\{1/F_{\omega}(z): |z|\le 1\right\}$ 
using \cref{lem:Lubich}(i).
This implies that $\zeta F_\omega(z)\neq 1$ for $|z|\leq 1$.  Note that $F_\omega(1)$ should be interpreted as the limit $\lim_{z\to 1}F_\omega(z)$, which shows that $\{0\} \notin \mathcal{S}_h^{\alpha}$, as expected.  
On the other hand, the condition that the F-LMM is strongly stable yields that $\delta(z)$ (cf. \eqref{eq:delta})  
is analytic and no zeros 
lie in a neighborhood of the closed unit disc $|z|\leq 1$ with the exception $z=1$.
Therefore, we know $F_{\omega}(z)=\frac{1}{(\delta(z))^{\alpha}}$ is also analytic for 
$z\in\Delta(R,\theta)$. This enables us to extend the result that $\zeta F_\omega(z)\neq 1$ 
from the unit disc $|z|\leq 1$ (with the exception $z=1$) to the lager region $z\in\Delta(R,\theta)$, hence 
verifies the validity of the condition \eqref{eq:stabcondition}.

Furthermore, by means of the strong stability, we know 
$F_\omega(z)$ has the  factorization representation
\begin{gather}\label{eq:Fexpan2}
F_\omega(z)=(\delta(z))^{-\alpha} = (1-z)^{-\alpha}F_{1}(z),
\end{gather}
where $F_{1}(z)$ is holomorphic at $z=1$ and $F_{1}(1)\neq0$ (see, e.g., \cite{Hairer}). 
From this expression, we readily see 
$
F_\omega(z)\sim C(1-z)^{-\alpha}~\text{ as } z\to 1,
$
with the constant $C=F_{1}(1) \neq 0$. Furthermore, using the constancy assumption that $F_\omega(z)=\delta(z)^{-\alpha}$ and 
$\delta\left(\mathrm{e}^{-h}\right)/h=1+O\left(h^{p}\right)$ as $h\to 0$ with $p\geq1$ in \eqref{eq:assum}, 
we find that $C=F_{1}(1)=1$.
It follows from \cref{lem:generating} that
\begin{equation} \label{eq:scalarexpan}
\begin{split}
F_y(z)&\sim y_0 \left(1+\frac{z}{(1-z)-\zeta (1-z)^{1-\alpha}}\right)\\ 
&= y_0 \frac{1}{(1-z)^{1-\alpha}} \left((1-z)^{1-\alpha}+\frac{z}{(1-z)^{\alpha}-\zeta }\right)\\
&\sim -\frac{y_0}{\zeta }\cdot \frac{1}{(1-z)^{1-\alpha}} \quad \text{as}\quad z\to 1,\\
\end{split}
\end{equation}
which leads to the desired result for Case I:
\begin{equation} \label{eq:scalardecay}
\begin{split}
y_n\sim -\frac{y_0}{\lambda  \Gamma(1-\alpha)}\cdot h^{-\alpha}n^{-\alpha}=O(t_{n}^{-\alpha}), \quad 
n\,\to \infty\,.
\end{split}
\end{equation}

{\bf Case II:} $\mathcal{L}$1 scheme for the scalar test equation. 
From the generating function \eqref{eq:generatingfunL1} for $\mathcal{L}$1 scheme, we should first derive 
the asymptotic behavior of the polylogarithm function $\mathrm{Li}_{\gamma}(z)$ as $ z\to 1$.
We know $\mathrm{Li}_{\gamma}(z)$ is well defined for $|z|<1$, can be analytically extended to 
$\mathbb{C}\setminus [1, \infty)$, and has the singular expansion \cite[Theorem VI.7]{FS09}:
\begin{equation} \label{eq:Liexpan}
\begin{split}
\mathrm{Li}_{\gamma}(z)\sim \Gamma(1-\gamma)w^{\gamma-1}+\sum_{j=0}^{\infty}\frac{(-1)^j}{j!}\zeta(\alpha-j)w^{j} 
\end{split}
\end{equation}
for all $\gamma\notin\{1,2,...\}$, with $\zeta(s)=\sum_{n=1}^{\infty}{1}/{n^s}$ (the Riemann zeta function) and 
$w=\sum_{\ell=1}^{\infty} {(1-z)^\ell}/{\ell}$. 
In particular, we know that the main asymptotic term of $\mathrm{Li}_{\gamma}(z)$ for $\gamma\in(0,1)$  
is given by $\mathrm{Li}_{\gamma}(z)\sim \Gamma(1-\gamma)(1-z)^{\gamma-1}$ as $z\to 1$.
Hence, the generating function $F_\mu(z)$ of the $\mathcal{L}$1 scheme has the asymptotic behavior 
\begin{gather}\label{eq:L1expan}
F_\mu(z)=\frac{1}{\Gamma(2-\alpha)} \frac{(1-z)^2}{z} \mathrm{Li}_{\alpha-1}(z) \sim (1-z)^{\alpha} \quad \mbox{as} ~~z\to 1\,,
\end{gather}
which implies the expansion $F_\omega(z) \sim (1-z)^{-\alpha}$ as $z\to 1$
for the $\mathcal{L}$1 scheme.

Using \cref{lem:generating}, we should now check the analytic properties of $F_\mu(z)$ on $\Delta(R, \theta)$.
For this, we write $F_\mu(z)=(1-z)^2 F_2(z)$, with $F_2(z)={\mathrm{Li}_{\alpha-1}(z)}/{(z\Gamma(2-\alpha))}$. 
We can easily see that  $\lim_{z\to 0}F_2(z)= \lim_{z\to 0} \mathrm{Li}'_{\alpha-1}(z)/\Gamma(2-\alpha)=1/\Gamma(2-\alpha)$. 
Hence, $z=0$ is a removable singularity of $F_2(z)$. On the other hand, we can see that $\mathrm{Li}_{\alpha-1}(z)\neq0$ for $z\in\Delta(R, \theta)\setminus \{0\}$.
Therefore we can redefine $F_2(0)=1/\Gamma(2-\alpha)$ to get 
$\tilde{F}_{\mu}(z)=F_\mu(z)$ 
for $z\in\Delta(R, \theta)\setminus \{0\}$, and $\tilde{F}_{\mu}(z)=1/\Gamma(2-\alpha)$ at $z=0$.
Now function $\tilde{F}_{\mu}(z)$ is analytic and has no zeros in $\Delta(R, \theta)$. Hence,
$\tilde{F}_{\omega}(z)=(\tilde{F}_{\mu}(z))^{-1}$ is also analytic
in $\Delta(R, \theta)$.
Then the desired results follows from \cref{lem:generating}. This complete the proof of Case II. 

{\bf Case III:} F-LMM and $\mathcal{L}$1 scheme for the general vector-valued system.
We now extend the previous proofs of Cases I and II for the scalar test equation to the vector-valued system. 
First for the F-LMM, we can multiply both sides of equation \eqref{eq:FLMMH2} with $z^{n}$ and 
then sum over $n\ge 0$ to obtain
\begin{equation*} 
\begin{split}
F_y(z)= (1-z)^{-1}y_0+h^{\alpha} A \left( F_\omega (z) F_{y} (z)- F_\omega (z) y_0 \right)\,.
\end{split}
\end{equation*}
This formula implies the representation of the solution to the F-LMM:
 \begin{equation} \label{eq:soluv}
\begin{split}
F_y(z)&= \left( I- h^{\alpha}F_{\omega}(z) A \right)^{-1} \left( (1-z)^{-1}I- h^{\alpha}F_{\omega} (z) A \right) y_0 \\
&= \left(I+\frac{z}{1-z }(I- h^{\alpha}F_{\omega}(z) A)^{-1} \right) y_0.
\end{split}
\end{equation}

For the $\mathcal{L}$1 scheme, we can define that $F_{\omega}(z)=(F_{\mu}(z))^{-1}$ 
by the convolution inverse and redefining $F_{\mu}(0)=1/\Gamma(2-\alpha)$ at the removable singularity point 
$z=0$ (see \eqref{eq:generatingfunL1} for $F_{\mu}(z)$). 
Hence, we see the formula \eqref{eq:soluv} holds for both F-LMM and $\mathcal{L}$1 scheme.

We know from \cref{thm:stability} that the inverse of $I- h^{\alpha}F_{\omega}(z) A$ 
exists and $\left( h^{\alpha}F_{\omega}(z) \right)^{-1}$ is not an eigenvalues of $A$ on 
$|z|\leq 1$ for the strong $A$-stable F-LMM or $\mathcal{L}$1 scheme. 
Hence, $F_3(z):=(I- h^{\alpha}F_{\omega}(z) A)^{-1}$ is analytic on 
$|z|\leq 1$, with exception $z=1$. We rewrite $F_3(z)= h^{-\alpha}F_{\mu}(z)\cdot (h^{-\alpha}F_{\mu}(z) I- A)^{-1}$. 
Following exactly the same argument as for the scalar case, the strong stability condition for the F-LMM 
and the structure of the generating function for the $\mathcal{L}$1 scheme 
enable us to verify the analyticity of $F_3(z)$ on $z\in \Delta(R,\theta)$, implying 
the analyticity of $F_y(z)$ for $z\in \Delta(R,\theta)$. 

Note that for both F-LMM and $\mathcal{L}$1 scheme we still have the asymptotical expansion 
$F_\omega(z)\sim (1-z)^{-\alpha}$ as $z\to 1$. Hence, we can derive 
by the existence of $A^{-1}$ and the fact that $\lambda_{A}\in\Lambda_{\alpha}^{s}$ 
\begin{equation} \label{eq:vectorexpan}
\begin{split}
F_y(z)&\sim  \left(I+\frac{z}{1-z } \left( I- \frac{h^{\alpha}}{(1-z)^{\alpha}} A \right)^{-1} \right) y_0\\
&= \left(I+\frac{z }{(1-z)^{1-\alpha} } \frac{1}{h^{\alpha}}  \left( \frac{(1-z)^{\alpha}}{h^{\alpha}} I-  A \right)^{-1} \right) y_0\\
&= \frac{1 }{(1-z)^{1-\alpha} } \left((1-z)^{1-\alpha} I+\frac{z}{h^{\alpha}}  \left( \frac{(1-z)^{\alpha}}{h^{\alpha}} I-  A \right)^{-1} \right) y_0\\
&\sim  \frac{1}{(1-z)^{1-\alpha}} \left(-\frac{1}{h^{\alpha}} A^{-1} \right) y_0 \quad \text{as}\quad z\to 1\,.
\end{split}
\end{equation}
Then we can readily see from \cref{lem:generating} the following result as we expect for Case III:
\begin{equation} \label{eq:vectordecay}
\begin{split}
y_n\sim -\frac{1}{  \Gamma(1-\alpha)}A^{-1}y_{0}\cdot h^{-\alpha}n^{-\alpha}, \quad \text{i.e.},~  \|y_n\|\sim O( t_{n}^{-\alpha}) \quad \text{as}\quad n\to \infty.
\end{split}
\end{equation}
\end{proof}

\section{Mittag-Leffler stability with small perturbations}
\label{sec:Perturbations}
The Mittag-Leffler stability of trivial solutions to F-ODEs with small perturbations 
(cf.\,\cref{lem:Cong}) was proved in \cite{Cong19} 
by the fixed-point technique of the Lyapunov-Perron operator. 
We can not see a possibility to apply the analysis in \cite{Cong19} 
for the study of the Mittag-Leffler stability of numerical solutions.
To motivate our subsequent analysis of numerical Mittag-Leffler stability,
we will recall the derivation of the variation of constants formula for F-ODEs 
by the Laplacian transformation and demonstrate that the fractional resolvent operators 
can be represented by the Mittag-Leffler functions.

Let $\widehat{f}=\mathcal{L}(f)(z)$ be the Laplace transform. 
For the Caputo fractional derivative with order $\alpha\in(0, 1)$, we have the identity  
$\widehat{_{~0}\mathcal{D}_{t}^{\alpha}y}(z)=z^{\alpha} \widehat{y}(z)-z^{\alpha-1}y(0)$ \cite{Diethelm10, Pod}.
Applying Laplace transform to the model \eqref{model} leads to that 
$z^{\alpha} \widehat{y}(z)=A\widehat{y}(z) +\widehat{f}(z)+z^{\alpha-1}y(0),$
that is, 
$
\widehat{y}(z)=(z^\alpha I-A)^{-1} \left( \widehat{f}(z)+z^{\alpha-1}y(0) \right).
$
By inverse Laplace transform and the convolution rule, the solution $y(t)$ can be represented by 
\begin{equation} \label{eq:SoluF}
\begin{split}
y(t)=R_{\alpha,1}(t)y_0+\int_{0}^{t} R_{\alpha,\alpha}(t-s)f(s, y(s))ds,
\end{split}
\end{equation}
where the fractional resolvent operator family $R_{\alpha,\beta}(t)$  generated by $A$ is defined by 
\begin{equation} \label{eq:FGoperators}
\begin{split}
R_{\alpha,\beta}(t):=\mathcal{L}^{-1} \left( z^{\alpha-\beta} (z^\alpha I-A)^{-1} \right) (t)=\frac{1}{2\pi i}\int_{\mathcal{C}} e^{zt}z^{\alpha-\beta}(z^\alpha I-A)^{-1}dz,
\end{split}
\end{equation}
$(\beta=1, \alpha)$ and a common choice for the integral path $\mathcal{C}$ is a line linking $c-i\infty$ and $c+i\infty$, where $c>0$ is a given positive constant.
In order to get the long-time decay rate of $\|y(t)\|$, the key is to derive sharp estimate of 
the operators $R_{\alpha,1}(t)$ and $R_{\alpha,\alpha}(t)$. 
In fact, $R_{\alpha,\beta}(t)$ can be expressed by the classical Mittag-Leffler function $E_{\alpha, \beta}$ and has the polynomial decay rate given by
 \begin{equation} \label{eq:FGisMF}
\begin{split}
R_{\alpha,1}(t) =& E_{\alpha}(t^{\alpha}A), \quad R_{\alpha,\alpha}(t) = t^{\alpha-1} E_{\alpha, \alpha}( t^{\alpha} A) ~\text{for} ~t>0,\\
\|R_{\alpha,1}(t)\|=& \|E_{\alpha}(t^{\alpha}A)\|=O( t^{-\alpha}), \quad \|R_{\alpha,\alpha}(t)\|= \|t^{\alpha-1} E_{\alpha, \alpha}( t^{\alpha} A)\|=O( t^{-\alpha-1})~ \text{as}~ t\to \infty,
\end{split}
\end{equation}
which lead to the variation of constants formula
\begin{equation} \label{eq:VCF}
\begin{split}
y(t)=E_{\alpha}(t^{\alpha}A) y_0 +\int_{0}^{t} (t-s)^{\alpha-1} E_{\alpha, \alpha}( (t-s)^{\alpha} A) f(s, y(s))ds.
\end{split}
\end{equation}

Although \eqref{eq:VCF} is well known, we are not aware of the appropriate literature to provide a derivation of \eqref{eq:VCF}  from \eqref{eq:FGoperators}. 
The standard approach is to use the Laplace transform and its inverse transform \cite{Pod}, but we emphasize here that we get \eqref{eq:VCF} directly from \eqref{eq:FGoperators}. This is because the inverse of discrete Laplace transform, or corresponding generating functions, that we derive from numerical methods are not known and thus require a fine estimate from the corresponding integral expression formula.
 Note that it is critical here to keep the optimal long-time decay rate of the resolvent operator.
The equations \eqref{eq:FGisMF} can be obtained in the following way:
\begin{equation} \label{eq:InverTM1}
\begin{split} 
R_{\alpha, \beta}(t) &=\frac{1}{2\pi i}\int_{\mathcal{C}} e^{tz} z^{\alpha-\beta}(z^{\alpha}I-A)^{-1} dz
=\frac{1}{2\pi i}\int_{\mathcal{C}} e^{tz} z^{\alpha-\beta} z^{-\alpha} (I-z^{-\alpha}A)^{-1} dz\\
&=\frac{1}{2\pi i}\int_{\mathcal{C}} e^{tz} z^{\alpha-\beta} z^{-\alpha} \sum_{k=0}^{\infty}(z^{-\alpha} A)^{k} dz 
=\sum_{k=0}^{\infty} \left( \frac{1}{2\pi i}\int_{\mathcal{C}} e^{tz} z^{\alpha-\beta} z^{-\alpha} z^{-k\alpha}  dz\right) A^{k},
\end{split}
\end{equation}
where we have chosen {the positive constant $c$ in}
the contour $\mathcal{C}$ large enough such that $\|z^{-\alpha} A\|<1$ and used the identity $(I-M)^{-1}=\sum_{k=0}^{\infty}M^{k}$ for any matrix $\|M\|<1$.
Notice that the integrand in the above last integral has only the singularity $z=0$, so the contour $\mathcal{C}$ can be deformed to $\Gamma_{(r, \theta)}$, 
where  
\begin{equation} \label{eq:counter}
\begin{split}
\Gamma_{(r, \theta)}:=\{z\in \mathbb{C}: |z|=r, \arg(z)\leq \theta\} \cup \{z\in \mathbb{C}: z=\rho e^{\pm i\theta}, \rho\geq r \}, \quad r>0,~ \frac{\pi}{2}<\theta\leq \pi.
\end{split}
\end{equation}
Applying the classical contour integral representation for the reciprocal Gamma function \cite[Page 14]{Pod}
\begin{equation} \label{eq:Gammafun}
\begin{split} 
\frac{1}{\Gamma(z)}&=\frac{1}{2\pi i} \int_{\Gamma(r, \theta)} e^{u} u^{-z} du 
=\frac{1}{2\pi\alpha i} \int_{\Gamma(r, \theta)} \exp{(\xi^{\frac{1}{\alpha}})} \xi^{(1-z-\alpha)/\alpha} d\xi,
\end{split}
\end{equation}
we get
\begin{equation} \label{eq:InverTM2}
\begin{split} 
R_{\alpha, \beta}(t) 
=\sum_{k=0}^{\infty} \left( \frac{1}{2\pi i}\int_{\Gamma(r, \theta)} e^{tz} z^{\alpha-\beta} z^{-\alpha} z^{-k\alpha}  dz\right) A^{k}
= t^{\beta-1}\sum_{k=0}^{\infty} \frac{(t^{\alpha}A)^{k}}{\Gamma(k\alpha+\beta)}= t^{\beta-1} E_{\alpha, \beta}( t^{\alpha} A).
\end{split}
\end{equation}

The equation \eqref{eq:InverTM2} implies readily that $R_{\alpha, 1}(t) = E_{\alpha}(t^{\alpha}A)$ and $R_{\alpha, \alpha} = t^{\alpha-1} E_{\alpha, \alpha}( t^{\alpha} A)$.  
We know that $\|E_{\alpha}(t^{\alpha}A)\|=O( t^{-\alpha})$ and $\|t^{\alpha-1} E_{\alpha, \alpha}( t^{\alpha} A)\|=O( t^{-\alpha-1})$ as $t\to \infty$ for $\lambda_A\in\Lambda_\alpha^s$ \cite{Cong18}.
Combining these estimates with the variation of constants formula \eqref{eq:VCF}, 
we can get the long-time optimal decay rate of the solution $\|y(t)\|=O( t^{-\alpha})$ under some smallness assumption of $f(t, y)$, that is, the Mittag-Leffler stability \cite{Cong18,Cong19}.

We like to point out that there is another effective way to estimate the continuous resolvent operator $R_{\alpha,\beta}(t)$ in \eqref{eq:FGoperators} or its discrete version, namely, to apply the following standard resolvent estimation formula 
\begin{equation} \label{eq:resestimate}
\begin{split}
\|(z^{\alpha}I -A)^{-1}\|\leq C_{\phi}|z^{\alpha}|^{-1},\quad  \forall z^{\alpha}\in \Sigma_{\phi}, \phi\in(0, \pi),
\end{split}
\end{equation}
where $\Sigma_{\theta}:=\{ z\in \mathbb{C}\setminus \{0\}: \arg(z)\leq \theta\}$.
Though this estimate is very effective in numerical analysis on the finite interval $[0, T]$ for fixed $T>0$ \cite{Jin16,Jin17}, it is not accurate enough to derive the long-time optimal decay rate, 
in particular, it will not enable us to achieve the discrete Mittag-Leffler stability with our 
desired decay rate for the numerical solutions. For a more detailed explanation, see \cref{sec:appendixB}.

The above idea motivates us that we may also define the corresponding discrete fractional resolvent operator family corresponding to $R_{\alpha, \beta}(t)$ for the numerical solutions based on the generating function, and then derive an accurate estimate of the discrete operators and obtain the Mittag-Leffler stability of the numerical solutions. In fact, this generating function approach has been widely used in various numerical analysis for time fractional differential equations \cite{Jin16,Jin17,LiWang19}.

\subsection{Discrete fractional resolvent family for F-LMMs and their decay rate} 
\label{sec:perturFLMMs}
We now study the numerical Mittag-Leffler stability of F-LMMs for the F-ODE model \eqref{model} with small perturbations.
Multiplying the equation \eqref{eq:FLMMs} with $z^{n}$ and summing the resulting equations over $n$ 
from $n=0$ to $\infty$, we obtain the generating function of F-LMMs: 
\begin{equation*} 
\begin{split}
F_y(z)=(1-z)^{-1}y_{0}+ h^{\alpha} A \left( F_\omega (z) F_{y} (z)-F_\omega (z) y_0 \right)+h^{\alpha}  \left( F_\omega (z) F_{f} (z)- F_\omega (z) f_0 \right),
\end{split}
\end{equation*}
where $f_n=f(t_n, y_n)$ and $F_{f} (z)=\sum_{n=0}^{\infty}f_{n}z^{n}$.
The above formula admits the solution representation 
 \begin{equation} \label{eq:pertursoluv}
\begin{split}
F_y(z)=& \left( I- h^{\alpha}F_{\omega}(z) A \right)^{-1} y_0 \left( (1-z)^{-1}I- h^{\alpha}F_{\omega} (z) A \right)
+h^{\alpha}  F_\omega (z)  \left(F_{f} (z)-f_0 \delta_{d} \right) \left( I- h^{\alpha}F_{\omega}(z) A \right)^{-1}\\
=& \left(I+\frac{z}{1-z }(I- h^{\alpha}F_{\omega}(z) A)^{-1} \right)y_0+ h^{\alpha}  F_\omega (z) \left( I- h^{\alpha}F_{\omega}(z) A \right)^{-1} \left(F_{f} (z)-f_0 \delta_{d} \right) \\
:=&F_{d}(z) y_0 +F_{D}(z)\left(F_{f} (z)-f_0 \delta_{d} \right),
\end{split}
\end{equation}
where $F_{d}(z)$ and $F_{D}(z)$ are given by 
 \begin{equation} \label{eq:FdD}
\begin{split}
F_{d}(z):&=I+\frac{z}{1-z }(I- h^{\alpha}F_{\omega}(z) A)^{-1} =\sum_{n=0}^{\infty} d_{n}z^{n},\\
F_{D}(z):&=  h^{\alpha} F_\omega (z) \left( I- h^{\alpha}F_{\omega}(z) A \right)^{-1}=\sum_{n=0}^{\infty} D_{n}z^{n}.\\
\end{split}
\end{equation}

In order to derive the discrete constant variation formula, we first note that 
 \begin{equation} \label{eq:discreConvo}
\begin{split}
F_{D}(z)\left(F_{f} (z)-f_0 \delta_{d} \right)&= h^{\alpha} F_\omega (z) \left( I- h^{\alpha}F_{\omega}(z) A \right)^{-1} \left(F_{f} (z)-f_0 \delta_{d} \right)\\
& = \left( \sum_{n=0}^{\infty} D_{n}z^{n} \right) \left(  \sum_{n=0}^{\infty} \tilde{f}_{n}z^{n}  \right)
= \sum_{n=0}^{\infty}  \left(\sum_{k=0}^{n} D_{n-k} \tilde{f}_{k}  \right) z^{n},  \\
\end{split}
\end{equation}
with $\tilde{f}=(0, f_1, f_2,...)$.  
Using this relation and comparing the coefficients of $z^n$ in \eqref{eq:pertursoluv}, we can get  
 \begin{equation} \label{eq:dVCF}
\begin{split}
y_{n}=d_{n}y_{0}+  \sum_{k=0}^{n} D_{n-k} \tilde{f}_{k}= d_{n}y_{0}+ \sum_{k=1}^{n} D_{n-k} f_{k}, \quad n\geq 1. \\
\end{split}
\end{equation}
Using the Cauchy formula, we can compute the coefficients $d_{n}$ and  $D_{n}$ by
\begin{equation} \label{eq:dnestimate}
\begin{split}
d_{n}&=\frac{1}{2\pi i}\int_{|z|=\rho} \frac{1}{z^{n+1}} F_{d}(z) dz=\frac{1}{2\pi i}\int_{|z|=\rho}  \frac{ 1}{z^{n+1}} \left( I+\frac{z}{1-z }(I- h^{\alpha}F_{\omega}(z) A)^{-1} \right) dz\\
    &=\frac{1}{2\pi i}\int_{|z|=\rho}  \frac{ 1}{z^{n+1}} \frac{ 1}{1-z} \left( (1-z)I+zh^{-\alpha}F_{\mu}(z)  \left( h^{-\alpha}F_{\mu}(z) I-  A \right)^{-1} \right) dz,\\
\end{split}
\end{equation}
\begin{equation} \label{eq:Dnestimate}
\begin{split}
D_{n}&=\frac{1}{2\pi i}\int_{|z|=\rho} \frac{1}{z^{n+1}} F_{D}(z) dz=\frac{1}{2\pi i}\int_{|z|=\rho}  \frac{ 1}{z^{n+1}} \left(  h^{\alpha} F_\omega (z) \left( I- h^{\alpha}F_{\omega}(z) A \right)^{-1} \right)dz\\
    &=\frac{1}{2\pi i}\int_{|z|=\rho}  \frac{ 1}{z^{n+1}} \left( h^{-\alpha}F_{\mu}(z) I-  A \right)^{-1} dz,\\
\end{split}
\end{equation}
where $\rho$ is a small positive constant. 
The operators $d_{n}$ and  $D_{n}$ can be seen as the discrete approximation of 
the fractional resolvent family $R_{\alpha,1}(t)$ and $R_{\alpha,\alpha}(t)$ at $t=t_n$ 
(it is also called discrete fractional resolvent family \cite{Lizama17,Ponce20}).
We can see that these two operators serve as the discrete Mittag-Leffler functions.  
With the help of the discrete constant variation formula, a key step is 
to derive the sharp asymptotical behavior of these two discrete operators, 
similarly to the continuous case. 
By comparing the continuous variation of constant formula \eqref{eq:VCF} with 
the discrete version \eqref{eq:dVCF}, we can establish the following decay rates of $d_n$ and $D_n$.

\begin{lemma}\label{lem:dis-estimate}
If the F-LMMs are strongly $A$-stable, then there exists a constant $h_0>0$
such that for any $0<h<h_0$, the discrete operators $d_n$ and $D_n$ given in \eqref{eq:dnestimate} and \eqref{eq:Dnestimate} have the decay estimates
 \begin{equation} \label{eq:dis-estimate} 
\begin{split} 
\|d_{n}\|\leq C_{\alpha} t_{n}^{-\alpha},\quad
\|D_{n}\|\leq C_{\alpha} t_{n}^{-\alpha-1} ~\text{as}~ n\to \infty.
\end{split}
\end{equation}
\end{lemma}

To motivate our proof of this key lemma, we give some important remarks. 
We first know from the homogeneous case in  \cref{sec:homogenous} that $\|d_n\|\sim O(t_n^{-\alpha})$, 
and without the limitation $0<h<h_0$, it suffices to prove the second estimate in \eqref{eq:dis-estimate}.
There are two possible approaches: 
it is natural to first apply a singularity analysis of the generating function, similarly to the homogeneous case; 
the second one is to estimate the integral directly based on the expression \eqref{eq:Dnestimate} 
by means of the standard resolvent estimation formula \eqref{eq:resestimate}.

With the first approach, we need to show that $F_{D}(z)\sim C (1-z)^{-\beta}$ as $z\to 1$, with $\beta\neq \{0, -1, -2,...\}$; see \cref{lem:generating}. 
Taking F-BDF$1$ as a example,  we know $F_\omega (z)=(1-z)^{-\alpha}$ and can easily get 
\begin{equation} \label{eq:FDesti}
\begin{split}
\lim_{z\to 1} F_{D}(z)&=\lim_{z\to 1} h^{\alpha} (1-z)^{-\alpha} \left( I- h^{\alpha} (1-z)^{-\alpha} A \right)^{-1}
=\lim_{z\to 1} h^{\alpha} \left( (1-z)^{\alpha} - h^{\alpha} A \right)^{-1}=- A ^{-1}.
\end{split}
\end{equation}
Hence, $F_{D}(z)\sim - A ^{-1}(1-z)^{0}$ as $z\to 1$, which has a degenerate index $\beta=0$, so 
\cref{lem:generating} cannot be used.

We now recall the main steps of the second approach. 
By changing the variable $z=e^{-\xi h}$, one gets
\begin{equation} \label{eq:dcvf3}
\begin{split} 
D_{n}=\frac{1}{2\pi i} \int_{\Gamma^{h}} e^{z t_{n}}  \left(h^{-\alpha}  F_{\mu} (e^{-zh}) I- A \right)^{-1} dz,\\
\end{split}
\end{equation}
where $\Gamma^{h}$ is given by $\Gamma^{h}:=\{ z=-\ln(\rho)/h+ i y: y\in \mathbb{R}~ \hbox{and} ~ |y|\le \pi/h\}.$
The contour $\Gamma^{h}$ can be deformed to $\Gamma_{(r, \theta)}$ in \eqref{eq:counter}.
A key step is to apply the standard resolvent estimate \eqref{eq:resestimate} to control the term $\| \left(h^{-\alpha}  F_{\mu} (e^{-zh}) I- A \right)^{-1}\|$.
This technique of decomposition and estimation of contour integrals has become a powerful tool in the numerical analysis of time fractional equations \cite{Cuesta07, Jin16,Jin17,Lubich96}.
However, as we shall see, this approach is insufficient for us 
to derive the sharp long-time decay rate of $\|D_{n}\|$ as we expect. 
We refer to \cref{sec:appendixB} for more details why this standard method does not work by looking at F-BDF$1$ again as an example.

Through the careful analysis above, we now understand why we can not get the optimal decay rate in $\|D_{n} \|$, 
namely, the resolvent estimation in \eqref{eq:resestimate} is not accurate enough. 
In fact, the resolvent operator identity is usually an infinite series \cite[Page 37, eq. (5.6)]{Kato}, the so-called first Neumann series for the resolvent, whereas the inequality \eqref{eq:resestimate} only uses the first term of the infinite series. 
This reminds us to avoid using the resolvent estimate directly, therefore get rid of 
the reduction of the decay rate in the analysis.  
Therefore, we shall relate the discrete resolvent operator $D_{n}$ to the continuous one $R_{\alpha,\alpha}(t_{n})$ 
in some sense so that we can make use of the existing estimates of the continuous operator. This enables 
the discrete operator to exactly preserve the long-time optimal decay rate of the continuous version.

\begin{proof}[Proof of Lemma \ref{lem:dis-estimate}]
It follows from \eqref{eq:Dnestimate} that
\begin{equation} \label{eq:Dnest1}
\begin{split}
D_{n}& =\frac{1}{2\pi i}\int_{|z|=\rho}  \frac{ 1}{z^{n+1}} \left(  h^{\alpha} F_\omega (z) \left( I- h^{\alpha}F_{\omega}(z) A \right)^{-1} \right)dz\\
    &=\frac{1}{2\pi i}\int_{|z|=\rho}  \frac{ 1}{z^{n+1}} \left(  h^{\alpha} F_\omega (z) \sum_{k=0}^{\infty} (h^{\alpha} F_\omega (z)A)^{k} \right)dz \quad(\text{iff}~\|h^{\alpha} F_\omega (z)A\|<1)\\
    &=\frac{1}{2\pi i} \int_{\Gamma^{h}} e^{z t_{n}}  \left(  h^{\alpha} F_\omega (e^{-z h}) \sum_{k=0}^{\infty} (h^{\alpha} F_\omega (e^{-z h})A)^{k} \right)dz \quad(\text{by}~z=e^{-\xi h}), 
    \end{split}
\end{equation}
where in the second equality 
we have chosen the step size $0<h<h_0$ small enough such that $\|h^{\alpha} F_\omega (z)A\|<1$.
Next, we shall represent the term inside the bracket in the last integral in \eqref{eq:Dnest1}
as a power series with respect to $z$.
By the strong stability assumption, the function $F_\omega(z)$ has the factorization representation
$F_\omega(z)= (1-z)^{-\alpha}F_{1}(z)$,
where $F_{1}(z)$ is holomorphic at $z=1$ and $F_{1}(1)\neq0$.  
In fact, $F_\omega(z)$ is analytic in the closed neighborhood of the unit disk with the exception of an isolated 
singularity at $z=1$.
So the function $F_\omega (e^{-z h})$ is analytic with the exception of an isolated singularity at $z=0$.
Hence we can derive 
\begin{equation} \label{eq:Fomegaexpan}
\begin{split}
h^{\alpha} F_\omega (e^{-z h})&=h^{\alpha} (1-e^{-z h})^{-\alpha} F_{1}(e^{-z h})\\
&=h^{\alpha} \left( \sum_{k=1}^{\infty} (-1)^{k+1}\frac{(zh)^{k}}{k!} \right)^{-\alpha} \left(c_0+c_1z+O(z^2)\right) \\
&=c_0 z^{-\alpha}  \left( 1+\frac{\alpha h z}{2}+O(z^2) \right) \left(1+\frac{c_1}{c_0} z+O(z^2)\right)\\
&=c_0 z^{-\alpha}  \left( 1+d_{1}z +O(z^2) \right),\\
 \end{split}
\end{equation}
where we have used the fact that $(1+z)^{-\alpha}=1-\alpha z+O(z^{2})$, and 
$c_0=F_{1}(1)\neq0$, $d_1={\alpha h}/{2}+ {c_1}/{c_0}$.
Furthermore, using the constancy relationship that $F_\omega(z)=\delta(z)^{-\alpha}$ and 
$\delta\left(\mathrm{e}^{-h}\right)/h=1+O\left(h^{p}\right)$ as $h\to 0$ with $p\geq1$ in \eqref{eq:assum}, 
we find that $c_0=F_{1}(1)=1$.
Substituting the expansion \eqref{eq:Fomegaexpan} into \eqref{eq:Dnest1} yields 
\begin{equation} \label{eq:Dnest1x}
\begin{split}
D_{n}&=\frac{1}{2\pi i} \int_{\Gamma^{h}} e^{z t_{n}}  \left(  z^{-\alpha}  \left( 1+d_{1}z +O(z^2) \right)\cdot \sum_{k=0}^{\infty} \left( z^{-\alpha}  \left( 1+d_{1}z +O(z^2) \right)A \right)^{k} \right)dz\\
    &= \frac{1}{2\pi i} \sum_{k=0}^{\infty}  \int_{\Gamma^{h}} e^{z t_{n}} z^{-\alpha}  \left( 1+(k+1)d_{1}z+O(k^2z^2) \right)\cdot (z^{-\alpha} A)^{k} dz    
    :=I_{1}+I_{2}.
    \end{split}
\end{equation}
Next, we estimate the dominant term $I_1$ and the higher-order term $I_2$, one by one. 
We can rewrite $I_1$ as 
\begin{equation} \label{eq:Dnest2}
\begin{split} 
    I_{1}(t_n)&=\frac{1}{2\pi i} \sum_{k=0}^{\infty}  \int_{\Gamma_{(r,\theta)}} e^{z t_{n}}  \left(  z^{-\alpha} \left(  z^{-\alpha} A \right)^{k} \right)dz 
    = \sum_{k=0}^{\infty} \left( \frac{1}{2\pi i} \int_{\Gamma_{(r,\theta)}} e^{z t_{n}}    z^{-\alpha}  (z^{\alpha-\alpha})  z^{-k\alpha} dz \right)  A^{k}\\
    &= t_{n}^{\alpha-1} \sum_{k=0}^{\infty} \frac{(t_{n}^{\alpha} A)^{k}}{\Gamma(k\alpha+\alpha)}=t_{n}^{\alpha-1} E_{\alpha, \alpha}( t_{n}^{\alpha} A),
\end{split}
\end{equation}
by using the reciprocal Gamma function formula \eqref{eq:Gammafun}, 
where the integral path $\Gamma^{h}$ is deformed to $\Gamma_{(r,\theta)}$.  
From \eqref{eq:Dnest2}, we readily get 
$ \|I_{1}(t_n)\|=\|t_{n}^{\alpha-1} E_{\alpha, \alpha}( t_{n}^{\alpha} A)\|=O(t_{n}^{-\alpha-1})$.

For the term $I_2$, we can rewrite it as 
\begin{equation} \label{eq:Dnest3}
\begin{split} 
    I_{2}(t_n)&=\frac{d_1}{2\pi i}  \sum_{k=0}^{\infty}  \int_{\Gamma_{(r,\theta)}} e^{z t_{n}} z^{-\alpha} \left((k+1)z+O(k^2 z^2) \right)   (z^{-\alpha} A)^{k}  dz \\
    &=d_1  \sum_{k=0}^{\infty}  \left( (k+1) \frac{1}{2\pi i}  \int_{\Gamma_{(r,\theta)}} e^{z t_{n}}  z^{-\alpha} (z^{\alpha-(\alpha-1)}) z^{-k\alpha}  dz \right) A^{k} +R_{2} \\
 &=d_1  t_{n}^{\beta-1} \left(  \sum_{k=0}^{\infty} (k+1)  \frac{(t_{n}^{\alpha}A)^{k}}{\Gamma(k\alpha+\beta)}  \right)  +R_{2} \\
 &=d_1 t_{n}^{\beta-1} \left(  \sum_{k=0}^{\infty} \frac{ \Gamma(k+2)}{k! \Gamma(2)}  \frac{(t_{n}^{\alpha}A)^{k}}{\Gamma(k\alpha+\beta)}  \right) +R_{2} 
 = d_1 t_{n}^{\beta-1} E_{\alpha, \beta} ^{2} ( t_{n}^{\alpha} A)  +R_{2},
\end{split}
\end{equation}
with $\beta=\alpha-1$, where the
remainder term $R_2$ can be given by 
\[
R_{2}= \frac{1}{2\pi i}  \sum_{k=0}^{\infty}  \left( \int_{\Gamma_{(r,\theta)}} e^{z t_{n}}  z^{-\alpha}\cdot O(k^2 z^{2}) \cdot
 z^{-k\alpha}  dz \right) A^{k}.
\]
The function $E_{\alpha, \beta} ^{2} ( t_{n}^{\alpha} A)$ is the Prabhakar function with parameter $\gamma=2$; see \cref{sec:appendixA}. 
Using the identity that 
$E_{\alpha, \beta} ^{2} (z)=\left[ E_{\alpha, \beta-1}(z)+ (1-\beta+2\alpha)E_{\alpha, \beta}(z) \right] /{(2\alpha)} ]$,
we come to the approximation 
\begin{equation} 
\begin{split} 
d_1 t_{n}^{\beta-1} E_{\alpha, \beta} ^{2} ( t_{n}^{\alpha} A) & =\frac{d_1}{2\alpha}  t_{n}^{\alpha-2}  \left( E_{\alpha, \alpha-2}(t_{n}^{\alpha} A)+ (2+\alpha)E_{\alpha, \alpha-1}(t_{n}^{\alpha} A) \right)\\
&=d_1 t_{n}^{\alpha-2} \left( O(t_{n}^{-2\alpha})+O(t_{n}^{-2\alpha}) \right) =O(t_{n}^{-\alpha-2}).\\
\end{split}
\end{equation}
Similarly, we can find that the remainder term $R_2$ is a higher-order term. 
Therefore, we have shown $\|I_{2}\|=O(t_{n}^{-\alpha-2})$, hence concluded that 
$\|D_{n}\|=O(t_{n}^{-\alpha-1})$. 
\end{proof}

\subsection{ Discrete fractional resolvent family and Poisson transformation} 
The discrete fractional resolvent sequence for time fractional difference equations with step size $h=1$ has been an important tool to study the qualitative properties of the solutions to fractional difference equations, such as the $\ell_p$-maximum regularity and the existence and uniqueness \cite{Lizama17}. This concept has recently been extended to arbitrary step size $h>0$ in \cite{Ponce20}  and was used to construct numerical schemes for linear sub-diffusion equations. One main advantage of the $\alpha$-resolvent sequence is that it allows us to write numerical solutions in terms of discrete constant variation formulas, exactly like the continuous case given in \eqref{eq:SoluF}. 
At the same time, one can connect the continuous Caputo fractional derivative with the discrete $\alpha$-difference scheme through Poisson transformation, as well as the discrete fractional resolvent operator with the continuous one. 
In this way we can estimate the optimal decay rate of the discrete fractional resolvent operator by means of 
the properties of the continuous resolvent operator and the Poisson transformation.

\begin{definition}  \cite{Ponce20}
For any $\alpha\in(0,1)$ and sequence $v=(v_{0}, v_{1},...)$, 
the $\alpha$-fractional sum of $v$ with a stepsize $h>0$ is given by 
$\mathcal{J}_{h}^{\alpha}(v_{n}):=h^{\alpha} \sum_{j=0}^{n}k_{n-j}^{\alpha} v_{j}$ for $n\in\mathbb{N}_0$,
where the coefficients 
$k_{0}^{\alpha}=1$ and $k_{n}^{\alpha}=\frac{\Gamma(\alpha+n)}{\Gamma(\alpha)\Gamma(1+n)}$
for  $n\geq 1.$
Furthermore, the Caputo $\alpha$-fractional difference is given by 
\begin{equation} \label{fracdifference}
\begin{split} 
\widetilde{\mathcal{D}}_{h}^{\alpha}(v_{n}):= \mathcal{J}_{h}^{1-\alpha}\Big( \frac{v_{n}-v_{n-1}}{h} \Big) =\frac{1}{h^{\alpha}}\Big( \sum_{j=0}^{n}k_{n-j}^{1-\alpha} v_{j}-\sum_{j=0}^{n-1}k_{n-1-j}^{1-\alpha} v_{j}\Big), \quad n\geq1.
\end{split}
\end{equation}
\end{definition}

The Caputo $\alpha$-fractional difference operator gives the numerical scheme for the F-ODE model \eqref{model} as
 \begin{equation} \label{eq:Scheme2}
\begin{split} 
\widetilde{\mathcal{D}}_{h}^{\alpha}(y_{n})= Ay_n+f_n, \quad n\geq1.
\end{split}
\end{equation}
Note that the Caputo $\alpha$-fractional difference operator is not exactly equivalent to F-BDF$1$, but it 
differs only in the coefficients of its initial value, so it can be seen as a correction scheme of F-BDF$1$ with 
the initial value.

The $\alpha$-fractional sum and difference operators have several nice properties, of which 
the so-called Poisson transformation is one that was studied in \cite{Lizama17} 
and extended in \cite{Ponce20}.
For fixed $h>0$ and $n\in\mathbb{N}_0$, the discrete Poisson distribution  is given by 
\begin{equation} \label{Poisson}
\begin{split} 
\rho^{h}_{n}(t)=e^{-\frac{t}{h}} \left( \frac{t}{h} \right)^{n} \frac{1}{h n!}, \quad h>0, ~n\in\mathbb{N}_0.
\end{split}
\end{equation}
One can check that $\rho^{h}_{n}(t)\geq 0$, $\rho^{h}_{n}(t)=h^{-1} \rho_{n}(t/h)$, and 
$\int_{0}^{\infty} \rho^{h}_{n}(t)dt=1$ for all $n\in\mathbb{N}_0$, where $\rho_{n}(t)=e^{-t} {t^{n}}/{n!}$ is the standard Poisson transformation. The discrete Poisson distribution appeared early in \cite{Cuesta03}
and has been used to analyze qualitative properties of numerical solutions to 
integro-differential equations.
Poisson distribution is an effective tool to prove the following key discrete constant variation formula
for semi-linear F-ODEs.

\begin{lemma}\label{lem:dcvf}
(\cite{Ponce20})
There is a unique solution to the equation \eqref{eq:Scheme2} with initial value $y_0$, satisfying 
the following discrete constant variation formula
 \begin{equation} \label{eq:dcvf}
\begin{split} 
y_{n}=Q_{1}^{n}y_{0}+h\sum_{j=0}^{n}Q_{\alpha}^{n-j}f_{j}, \quad n\geq 1,
\end{split}
\end{equation}
where $Q_{\beta}^{n}$ is the discrete resolvent operator given by 
$Q_{\beta}^{n}=\int_{0}^{\infty} \rho^{h}_{n}(t)R_{\alpha,\beta}(t) dt$ ($\beta=1$ or $\alpha$) and $R_{\alpha,\beta}(t)=t^{\beta-1}E_{\alpha, \beta}(t^{\alpha}A)$ is the continuous fractional resolvent operator given in \eqref{eq:FGisMF}.
\end{lemma}

As we have seen in \cref{sec:perturFLMMs}, once the discrete constant variation formula is known, 
the next key step is to derive the sharp decay rate of the discrete fractional resolvent operators.
If $f\equiv0$, the contour integral representation for $Q_{1}^{n}$ can be derived by the standard discrete Laplace transform or a generating function approach \cite{Lubich85, Lubich86}. 
The long-time decay rate of $\|Q_{1}^{n}\|$ can be obtained by a singularity analysis using 
generating functions as that in \cref{sec:homogenous} or by a contour integral method 
\cite{Cuesta07}. However, neither method can lead to a desired optimal decay rate for $\|Q_{\alpha}^{n}\|$.

The nice relationship $Q_{\beta}^{n}=\int_{0}^{\infty} \rho^{h}_{n}(t)R_{\alpha,\beta}(t) dt$ presents
a completely new approach that can make full use of the properties of the discrete Poisson distribution $\rho_{n}^{h}(t)$ and the estimates of the continuous resolvent family $R_{\alpha,\beta}(t)$ (cf.\,\eqref{eq:FGisMF}).
The proof is very simple and elegant compared with the method in \cref{sec:perturFLMMs} developed for F-LMMs.  
However, this approach is effective only for this very specific scheme, for which the discrete and continuous 
Poisson transformations connecting the discrete resolvent operator to the continuous one happen to be known explicitly.

We now prove that $\|Q_{1}^{n}\|\leq {C}{/t_{n}^\alpha}$ and $\|Q_{\alpha}^{n}\|\leq {C}/{t_{n}^{1+\alpha}}$.
It follows readily from \cref{lem:dcvf} and the estimate $\left\| R_{\alpha,1}(t)\right\|\leq Ct^{-\alpha}$ that 
 \begin{equation} \label{eq:coeffestimate1}
\begin{split} 
\|Q_{1}^{n}\|&\leq \int_{0}^{\infty} \rho^{h}_{n}(t) \left\| R_{\alpha,1}(t)\right\| dt
\leq C   \int_{0}^{\infty} e^{-\frac{t}{h}} \left( \frac{t}{h} \right)^{n} \frac{1}{h n!} \frac{1}{t^{\alpha}} dt
=\frac{C}{n!}  \frac{1}{h^\alpha}  \int_{0}^{\infty} e^{-t} t^{n-\alpha} dt, \\
\end{split}
\end{equation}
Furthermore, we note the fact that 
${\Gamma(k+\alpha)k^{-\alpha}}/{\Gamma(k)} \to 1$ as $k\to\infty$ for any $\alpha>0$,
and hence can derive that 
 \begin{equation} \label{eq:coeffestimate1x}
\begin{split} 
\frac{1}{n!}  \frac{1}{h^\alpha}  \int_{0}^{\infty} e^{-t} t^{n-\alpha}  dt 
=\frac{1}{h^\alpha} \frac{\Gamma(n+1-\alpha) }{n\Gamma(n)} 
=\frac{1}{h^\alpha}\frac{n^{1-\alpha}}{n} \frac{\Gamma(n+1-\alpha) n^{\alpha-1}}{\Gamma(n)}
\leq \frac{C}{t_{n}^\alpha}~~\text{as}~n\to \infty,
\end{split}
\end{equation}
which leads to the estimate $\|Q_{1}^{n}\| \leq {C}/{t_{n}^\alpha}$.
Similarly, we can estimate $\|Q_{\alpha}^{n}\|$ by noting that 
$\left\| R_{\alpha,\alpha}(t)\right\|\leq Ct^{-\alpha-1}$, 
 \begin{equation} \label{eq:coeffestimate2}
\begin{split} 
\|Q_{\alpha}^{n} \|&\leq \int_{0}^{\infty} \rho^{h}_{n}(t) \left\| R_{\alpha,\alpha}(t)\right\| dt
\leq C   \int_{0}^{\infty} e^{-\frac{t}{h}} \left( \frac{t}{h} \right)^{n} \frac{1}{h n!} \frac{1}{t^{\alpha+1}} dt
=\frac{C}{n!}  \frac{1}{h^{\alpha+1}}  \int_{0}^{\infty} e^{-t} t^{n-\alpha-1}  dt.
\end{split}
\end{equation}
Now, applying the properties of the Gamma function again yields 
 \begin{equation} \label{eq:coeffestimate2x}
\begin{split} 
 \frac{1}{n!}  \frac{1}{h^{\alpha+1}}  \int_{0}^{\infty} e^{-t} t^{n-\alpha-1}  dt 
 =\frac{1}{h^{\alpha+1}} \frac{\Gamma(n-1+(1-\alpha)) }{n(n-1)\Gamma(n-1)} 
\leq \frac{C}{t_{n}}\frac{1}{t_{n-1}^{\alpha}} 
\leq \frac{C}{t_{n}^{\alpha+1}}~~\text{as}~n\to \infty, 
\end{split}
\end{equation}
which gives the desired estimate of $\|Q_{\alpha}^{n}\|$.

\subsection{Numerical Mittag-Leffler stability under perturbations} 
For deriving the desired optimal decay rate of numerical solutions, we still need the following lemma established 
in our early work \cite{Wang18}.
\begin{lemma}\label{lem:Rate}
Consider the Volterra difference equation 
\begin{equation*}
y_{n+1}= q_{n}+\sum\limits_{j=0}^{n}Q_{n-j}y_{j},~~~n\geq0.
\end{equation*}
If the coefficients satisfy 
 $q_{n}\sim \frac{c_{1}}{n^{\alpha}}, ~ Q_{n}\sim \frac{c_{2}}{n^{1+\alpha}}$  and $\rho:=\sum_{j=0}^{\infty}|Q_{j}|\leq \rho_0<1$
for some constants $c_{1},  c_{2}>0$ and $0<\alpha<1$, 
then the asymptotic estimate $y_{n}\sim \frac{c_{1}\left(1-\rho \right)^{-1}}{n^{\alpha}}$ is valid.
\end{lemma}

With all the preparations of this section up to now, we are ready to establish our main results.

\begin{theorem}\label{thm:perturbFODEs}
For the non-homogenous F-ODEs model \eqref{model}, we assume 
$\lambda_{A}\in \Lambda_{\alpha}^{s}$, and that $f$ is continuous,  $f(t, 0)=0$, and further satisfies 
\begin{equation} \label{eq:fcond2}
\begin{split}
\|f(t, x(t))-f(t, y(t))\|\leq L(t)\|x(t)-y(t)\|, \quad \forall\,t\geq 0, ~~x, y\in\mathbb{R}^{d}, 
\end{split}
\end{equation}
where $L: [0, \infty)\to \mathbb{R}_{+}$ is a continuous Lipschitz function. 
Letting $\mathcal{L}_0=\sup_{t\geq0} L(t)$ and parameter $D_0$ be defined as in \cref{tab:D0}, 
then there exists constant $h_0>0$ such that for any $0<h<h_0$, 
the trivial solution obtained by the strong A-stable F-LMM \eqref{eq:dVCF} or the $\alpha$-difference scheme \eqref{eq:dcvf} are numerically Mittag-Leffler stable, i.e., $\|y_{n}\|=O(t_n^{-\alpha})$ as $n\to \infty$, 
provided that the Lipschitz function $L(t)$ is small enough in the sense that 
\begin{equation} \label{eq:Lcond2}
\begin{split}
1-\|D_{0}\| \mathcal{L}_0 >0, \quad 
\frac{1}{ 1-\|D_{0}\| \mathcal{L}_0}  \Big( \lim_{n\to \infty} \sum_{k=0}^{n-1} \|D_{n-k}\| L(t_k) \Big)\leq \rho_0<1.
\end{split}
\end{equation}
\end{theorem}

\begin{proof}
For the $\alpha$-difference scheme \eqref{eq:dcvf}  (cf.\,\cref{lem:dcvf}), let $d_{n}=Q_{1}^{n}=\int_{0}^{\infty} \rho^{h}_{n}(t)R_{\alpha,1}(t)dt$ and $D_{n}=hQ_{\alpha}^{n}=h\int_{0}^{\infty} \rho^{h}_{n}(t)R_{\alpha,\alpha}(t)dt$.
Then both the F-LMM in \eqref{eq:dVCF} and the $\alpha$-difference scheme in \eqref{eq:dcvf} for solving the F-ODE model \eqref{model} with small perturbation can be written as a unified form:
\begin{equation} \label{eq:dVCF2}
\begin{split}
y_{n}= d_{n}y_{0}+ \sum_{k=1}^{n} D_{n-k} f_{k}, \quad n\geq 1, \\
\end{split}
\end{equation}
where the coefficients have the decays $\|d_{n}\|\leq C_\alpha t_{n}^{-\alpha}$ and  $\|D_{n}\|\leq C_\alpha t_{n}^{-\alpha-1}$ due to \cref{lem:dis-estimate} and the estimate \eqref{eq:coeffestimate1}-\eqref{eq:coeffestimate2x}, 
with $C_\alpha$ being independent of $t_n$.
It comes directly from the above equation that 
 \begin{equation} \label{eq:numsoluesti1}
\begin{split} 
\|y_{n}\|&\leq \|d_{n}\| \|y_{0}\|+\sum_{k=0}^{n} \|D_{n-k}\| \|f_{k}\|
\leq \|d_{n}\|  \|y_{0}\|+ \Big(\sum_{k=0}^{n-1} L(t_k) \|D_{n-k}\|  \|y_{k}\|  + \mathcal{L}_0 \|D_{0}\|  \|y_n\| \Big)
\end{split}
\end{equation}
for all $n\geq 1$, where the fact that 
$\|f_{k}\|=\|f(t, y_{k})-f(t, 0)\|
\leq \mathcal{L}_0 \|y_{k}\| $ is used.  This readily implies 
 \begin{equation} \label{eq:numsoluesti2}
\begin{split} 
\|y_{n}\| &\leq \frac{\|y_{0}\| } {1-\|D_{0}\| \mathcal{L}_0} \|d_{n}\| + \Big( \frac{1}{ 1-\|D_{0}\| \mathcal{L}_0}  \sum_{k=0}^{n-1} \|D_{n-k}\| L(t_k) \Big)  \|y_{k}\|.\\
\end{split}
\end{equation}
Now the Mittag-Leffler stability estimate $\|y_{n}\|=O(t_{n}^{-\alpha})$ as $t_n\to \infty$ follows directly 
from the decays of $\|d_{n}\|$ and $\|D_{n}\|$, the conditions \eqref{eq:Lcond2} and the asymptotic estimate 
in \cref{lem:Rate}.
\end{proof}

We give some remarks on the assumption on the smallness of 
the Lipschitz function in \eqref{eq:Lcond2}. 
First of all, if $L(t)$ is constant and $\mathcal{L}_0=L(t)$, then 
this assumption reduces to $\mathcal{L}_0< \frac{1}{\|D\|_0+S_0}$, where $S_0=\sum_{k=1}^{\infty} \|D_{k}\|$. Note that $S_0$ is finite due to the estimate $ \|D_{k}\|=O(t_{k}^{-\alpha-1})$.  
For five numerical methods studied in this work, 
their generating functions $F_{\omega}(z)$ and the values of $D_0$ are listed in \cref{tab:D0}.

On the other hand, by comparing the second condition in \eqref{eq:Lcond2} and the assumption \eqref{eq:Lcond}(i), we can find that $\lim_{n\to \infty} \sum_{k=0}^{n-1} \|D_{n-k}\| L(t_k)$ can be seen as a discrete version of $\int_{0}^{t} (t-s)^{\alpha-1}\|E_{\alpha, \alpha} ((t-s)^{\alpha}A)\| L(s)ds$ up to a constant. 
The condition \eqref{eq:Lcond2} for numerical methods is slightly stronger than the condition \eqref{eq:Lcond} 
for the continuous equation, but the results are also stronger, 
namely, it is asymptotically stable in the continuous case while it is Mittag-Leffler stable in the discrete case.

\begin{table}[h]
		\centering
		\caption{Numerical methods and their generating functions and parameters $D_0$}
		\label{tab:D0}
		\begin{tabular}{| l |c|c|c|}\hline
			Methods &\multicolumn{1}{c|}{$F_{\omega}(z)$}  &\multicolumn{1}{c|}{$F_{\omega}(0)$} &\multicolumn{1}{c|}{$D_{0}=h^{\alpha}F_{\omega}(0)(I-h^{\alpha}F_{\omega}(0)A)^{-1}$} \\\hline
			         
			F-BDF$1$&$ (1-z)^{-\alpha}$& 1 & $(h^{-\alpha}I-A)^{-1}$ \\\hline
			F-BDF$2$&$ (1-z)^{-\alpha} \left( \frac{3-z}{2}\right)^{-\alpha}$&$\left( \frac{2}{3}\right)^{\alpha}$& $\left( (\frac{2}{3})^{-\alpha}h^{-\alpha}I-A\right)^{-1}$ \\ \hline
			F-Adams$2$&$(1-z)^{-\alpha} \left(1-\frac{\alpha}{2}(1-z) \right)$&$1-\frac{\alpha}{2}$&$\left( (1-\frac{\alpha}{2})^{-\alpha}h^{-\alpha}I-A\right)^{-1}$ \\\hline
			$\mathcal{L}1$ method &$\frac{z} {(1-z)^2} \mathrm{Li}^{-1}_{\alpha-1}(z)$ & 1& $(h^{-\alpha}I-A)^{-1}$\\\hline
			$\alpha$-Difference &No explicit form&No explicit form& $\int_{0}^{\infty} e^{-\frac{t}{h}} t^{\alpha-1} E_{\alpha,\alpha}(t^\alpha A)dt$ \\\hline
	\end{tabular}		
\end{table}

We end this section with a general comment. Generating functions are an effective tool for the study of 
the long-term stability and convergence of numerical solutions to both integer and fractional evolution equations. However, various estimates based on the Gronwall-type inequality \cite{Kopteva19,Liao19} are mostly suitable only for numerical analysis over finite time, 
due to the common fact that the Gronwall-type  inequality often contains a growth factor of exponential or Mittag-Leffler functions, 
which is uncontrollable when time is not finite.  

\section{Applications and numerical examples} 
\label{sec:examples}
In this section, we present several representative examples to show the polynomial decay rate of numerical solutions obtained by numerical methods in \cref{tab:D0} for various time fractional F-ODEs, including 
the time fractional sub-diffusion equations, the fractional optical control system 
and the stable equilibrium points for nonlinear F-ODEs. 
In the concrete implementation of F-LMMs, it is very important to calculate the weight coefficients $\{\mu_{k}\}$ or $\{\omega_{k}\}$ effectively. It is generally not easy to quickly compute 
the coefficients of the fractional expansion of a rational polynomial function, but the Miller formula \cite[Theorem 4]{Garrappa15}
is very useful for the purpose. 
We have used the Miller formula to compute the coefficients of the F-BDF$k$ schemes and F-Adams$k$ schemes.

\subsection{Decay rate of F-ODEs} 
In this example, we consider the simple scalar F-ODE $_{~0}\mathcal{D}_{t}^{\alpha}y(t)=\lambda y$, 
with the eigenvalue $\lambda=1+(1+b)i, b\in\mathbb{R}$, which contains a positive real part, 
but the solution still decays polynomially. 
In order to test the numerical decay rate quantitatively, we introduce the index function 
\begin{equation} \label{eq:indexp}
\begin{split}
p_{\alpha}(t_{n})=-\frac{\ln(\|y_{n+m}\|/\|y_{n}\|)}{\ln(t_{n+m}/t_{n})}, 
\end{split}
\end{equation}
where $t_n>1$ and $m$ is a fixed integer. We take $m=5$ in all the following examples. The index $p_{\alpha}$ is a numerical observation of $\alpha$ given in $\|y_{n}\|=O(t_{n}^{-\alpha})$, which is independent of the initial value. 
In the simulation, we take the initial value $y(0)=5$ and the parameter $\lambda=1+(1+b)*i$ for $b=0.1, 0, -0.1$.

\begin{figure}[H]\label{fig01}
\begin{minipage}{1.00\linewidth}
\centering
\includegraphics[scale=0.4]{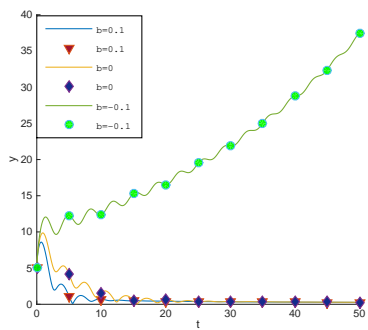}
\caption{Numerical solutions for $\alpha=0.5, h=0.1$ with $\lambda=1+(1+b)*i$ for $b=0.1, 0, -0.1$ computed by F-BDF1}
\end{minipage}
\end{figure}

For $\alpha=0.5$, we have $\Lambda_{0.5}^s=\left\{z\in\mathbb{C}\setminus \{0\}: |\arg(z)|> \frac{\pi}{4}\right\}$. Hence we can see $\lambda\in\Lambda_{0.5}^s$ for $b=0.1$, $\lambda\notin\Lambda_{0.5}^s$ for $b=-0.1$ and $\lambda\in\partial(\Lambda_{0.5}^s)$ for $b=0$. We observe from \cref{fig01} that  for both $b=0$ and $b=0.1$ the numerical solutions decay at a polynomial rate while for $b=-0.1$ the numerical solutions increase polynomially with time.  Other numerical methods, such as F-BDF$2$ and $\mathcal{L}1$ scheme, 
give almost the same results as this one, and therefore are not plotted here. This result shows the qualitative polynomial decay of solutions to F-ODEs, especially when the eigenvalue has a positive real part, 
which is very different from ODEs of integer order.

\begin{table}\label{table11}
\caption{Observed $p_\alpha$ computed by F-BDF$1$ and F-BDF$2$ (the data in the brackets) with $h=0.1$ and $b=10$ for Example 5.1} 
\begin{center}
\begin{tabular}{lll lll ll }
\hline $t_n$  & $\alpha=0.3$  & $\alpha=0.5$ & $\alpha=0.7$ & $\alpha=0.9$\\
\hline
$100   $ & 0.3009(0.3008)  & 0.5009(0.5008) & 0.7011(0.7009) & 0.9016(0.9011)\\
$200   $ & 0.3005(0.3005)  &0.5005(0.5004)  & 0.7005(0.7004)  & 0.9008(0.9006)\\
$300  $ & 0.3004(0.3004)  &0.5003(0.5003)  & 0.7003(0.7003)  & 0.9005(0.9004)\\
$400  $& 0.3004(0.3004)  & 0.5002(0.5002) & 0.7003(0.7002)  & 0.9004(0.9003)\\
$500 $& 0.3003(0.3003)  & 0.5002(0.5002) & 0.7002(0.7002)  & 0.9003(0.9002)\\
\hline
\end{tabular}
\end{center}
\end{table}

\begin{table}\label{table13}
\caption{ Observed $p_\alpha$ computed by $\mathcal{L}1$ scheme and F-Adams$2$ (the data in the brackets) with $h=0.1$ and $b=10$ for Example 5.1} 
\begin{center}
\begin{tabular}{lll lll ll }
\hline $t_n$  & $\alpha=0.3$  & $\alpha=0.5$ & $\alpha=0.7$ & $\alpha=0.9$\\
\hline
$100   $ &0.3008(0.3008)  &0.5008(0.5008) & 0.7009(0.7009)  & 0.9011(0.9011)\\
$200   $ &0.3005(0.3005)  &0.5004(0.5004)  & 0.7004(0.7004)  & 0.9006(0.9006)\\
$300  $ & 0.3004(0.3004)  &0.5003(0.5003)  &0.7003(0.7003)  & 0.9004(0.9004)\\
$400  $&0.3004(0.3004)  & 0.5002(0.5002) & 0.7002(0.7002)  & 0.9003(0.9003)\\
$500  $&0.3003(0.3003)  & 0.5002(0.5002) & 0.7002(0.7002)  & 0.9002(0.9002)\\
\hline
\end{tabular}
\end{center}
\end{table}

The energy method depends heavily on the special structure of coefficient $\{\mu_{j}\}_{j=0}^{\infty}$ \cite{Wang18}, 
and an additional requirement on step sizes is needed for the F-BDF$2$ schemes. 
The new results in this work show that the long-time polynomial decay of the numerical solutions 
is closely related to the numerical stability and there is no any step size restriction 
for the F-BDF$2$ schemes. This is also true for the F-Adams$2$ schemes.
In order to further quantitatively describe the decay rate of numerical solutions, we compute 
the observed index $p_\alpha$ with various choices of parameters
for F-BDF$1$, F-BDF$2$, F-Adams$2$ and $\mathcal{L}1$ schemes, as shown in
 \cref{table11} and \cref{table13} respectively, 
from which we can see that the numerical solutions decay clearly at the rate $O(t_{n}^{-\alpha})$, which 
is completely consistent with our theoretical prediction.

\subsection{Time fractional advection-diffusion equations} 
As the second example, we consider the time fractional advection diffusion problem: 
\begin{equation}\label{eq:admodel}
\begin{split}
_{~0}\mathcal{D}_{t}^{\alpha}u(x, t)+a\cdot \nabla u=D\Delta u, ~t>0, x\in\Omega=[0,1],\\
\end{split}
\end{equation}
with the initial value $u(x, 0)=u_{0}(x)$, the periodic boundary condition, and constant coefficients 
$a\in\mathbb{R}, D> 0$. 
If $u_{0}\in L^{2}(\Omega)$ and $u(x, t)=0$ for $x\in\partial\Omega$, it is known that 
the solution exhibits singularity near $t=0$, and 
$\|_{~0}\mathcal{D}_{t}^{\alpha} u(\cdot , t)\|_{L^2(\Omega)}\leq C_{\alpha} t^{-\alpha} \|u_{0}\|_{L^2(\Omega)}$, 
and more importantly, the solution decays in a polynomial rate,
 i.e., $\|u(\cdot , t)\|_{L^2(\Omega)}=O( t^{-\alpha})$ as $t\to +\infty$.

However, we are not aware of studies in the literature  
of the polynomial decay of solutions and characterizing their long-tail effect 
for time fractional advection diffusion equations from the numerical point of view.
{When $a=0$,} we know the eigenvalues of F-ODEs after spatial semi-discretizations of these equations 
are often negative real constants, therefore any $A_{0}$-stable numerical method \cite{Hairer},  i.e., the stable region contains the entire negative real half axis, are unconditionally stable. When $a\neq 0$, the F-ODEs contain eigenvalues with non-zero imaginary part, 
therefore $A(\frac{\pi}{2})$-stable numerical methods can overcome the restriction on step size due to stability. 
We did some simulations in \cite{LiWang19} for this example by $\mathcal{CM}$-preserving schemes and 
verified their {$A(\frac{\pi}{2})$-stability}.
We shall now confirm the long-term polynomial decay rate of the solution, namely, the Mittag-Leffler stability.

For the space discretization on a uniform grid $\{x_1, x_2, ...,x_N\}$ with grid points $x_j = j\delta x$ and mesh width $\delta x = 1/N$ under the periodic boundary condition $u(0, t)=u(1, t)$, we can apply the standard second-order central differences for the advection and diffusion terms in \eqref{eq:admodel} to get the semi-discrete system
\begin{equation}\label{eq:semiDSM}
\begin{split}
_{~0}\mathcal{D}_{t}^{\alpha}U(t)+ \frac{a}{2\delta x}BU= \frac{D}{{\delta x}^2} AU,\quad t>0,\\
\end{split}
\end{equation}
where $U(t)=(u_{1}, u_2,...,u_{N})^{T}$, $u_{0}=u_{N}, u_{N+1}=u_{1}$ and 
$$
B=\begin{pmatrix}
0 & 1 & 0& \cdots & -1 \\
-1 & 0 & 1& \cdots & 0 \\
 & \ddots & \ddots &\ddots & \\
0 & \cdots & -1& 0 & 1 \\
1 & 0 & \cdots& -1 & 0 \\
\end{pmatrix},\quad
A=\begin{pmatrix}
-2 & 1 & 0& \cdots & 1 \\
1 & -2 & 1& \cdots & 0 \\
 & \ddots & \ddots &\ddots & \\
0 & \cdots & 1& -2 & 1 \\
1 & 0 & \cdots& 1 & -2 \\
\end{pmatrix}
.$$

\begin{figure}[!ht]
\begin{center}
$\begin{array}{cc}
 \includegraphics[scale=0.4]{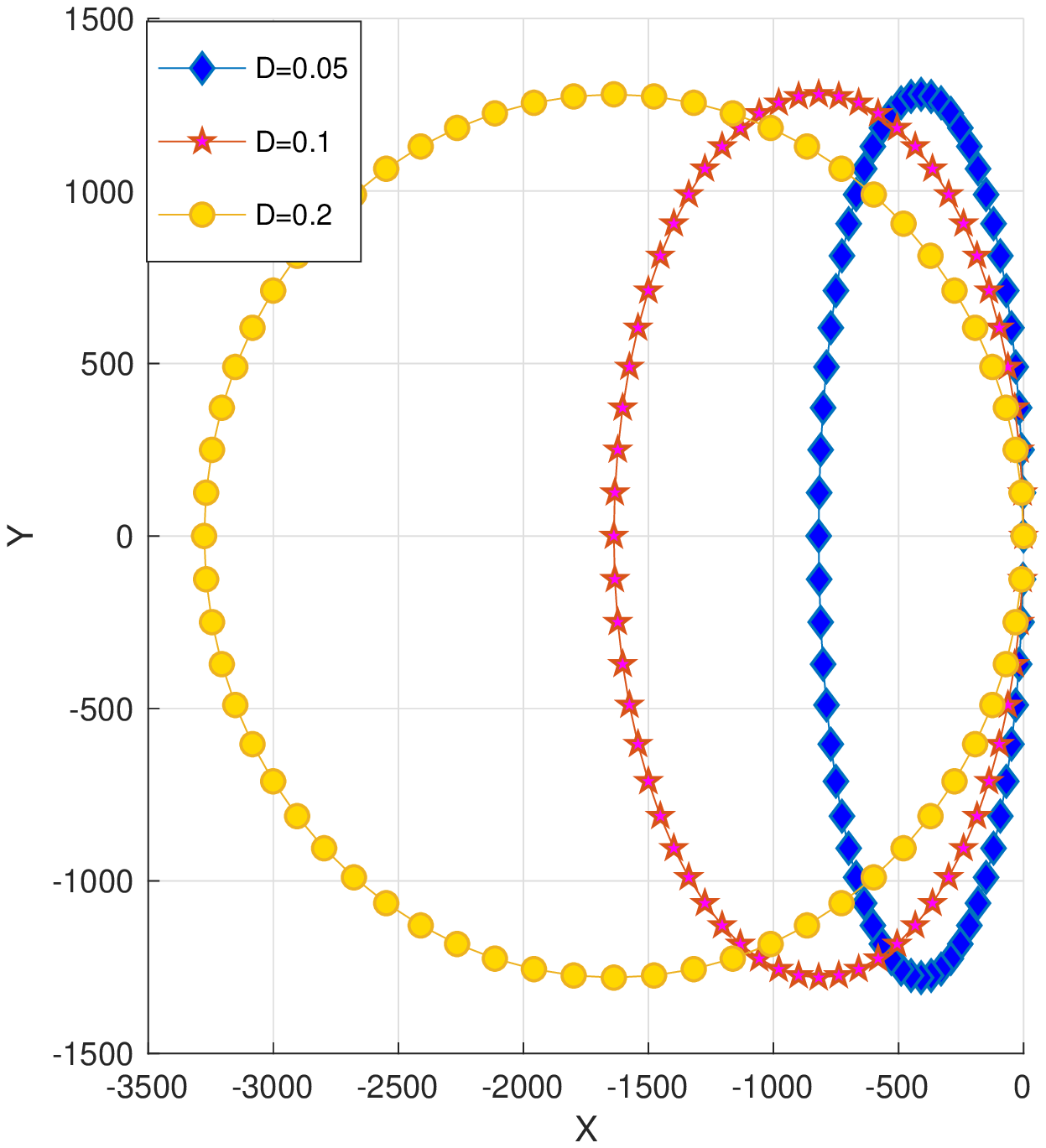} \qq
 \includegraphics[scale=0.4]{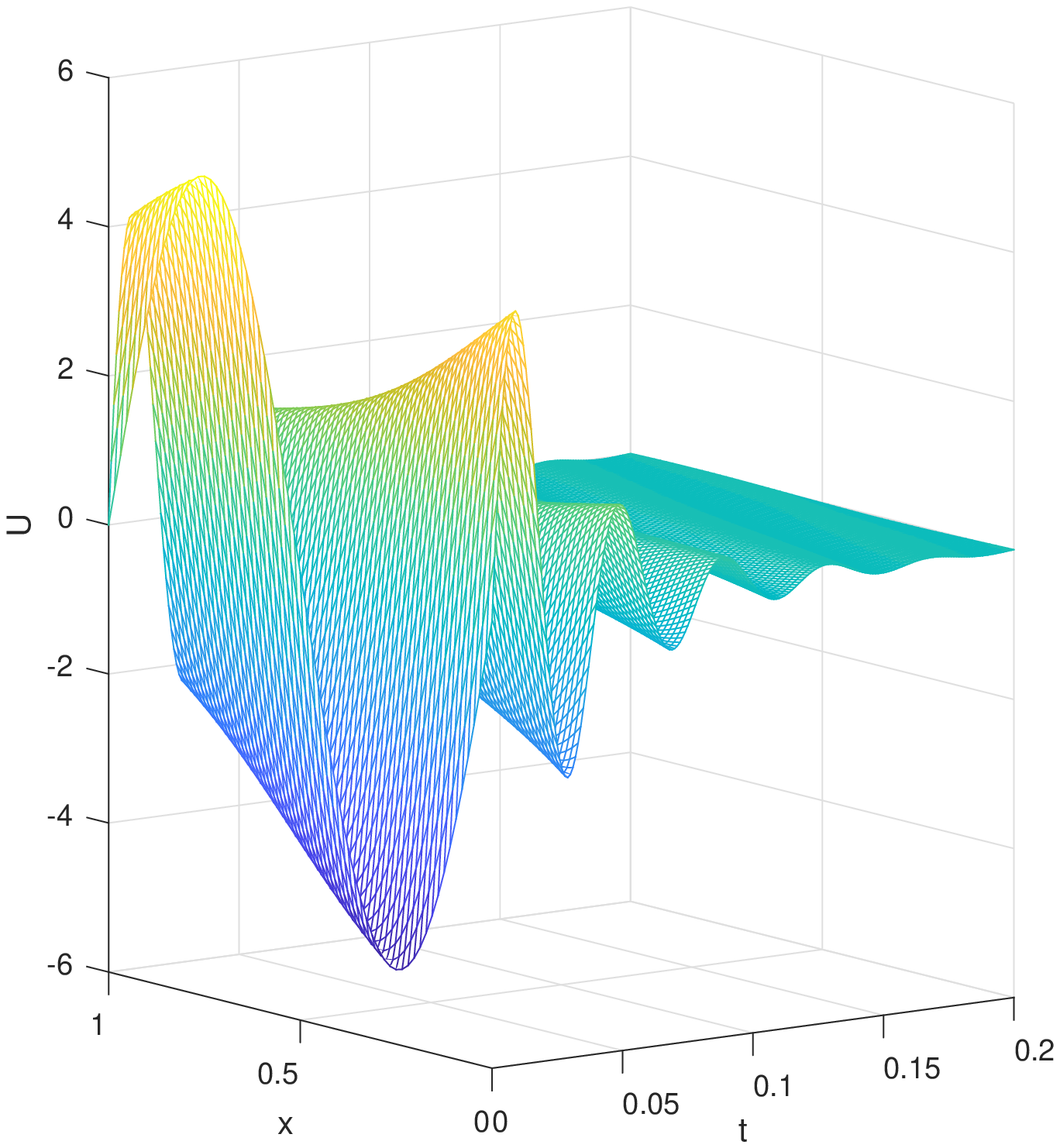}
 \end{array}$
\end{center}
\caption{Eigenvalues distributions with different parameters $D$ (Left); numerical solutions with $a=20, \delta x=1/64, \alpha=0.95, h=0.001$ and $u_0=5\sin(2\pi x)$ computed by $\mathcal{L}1$ scheme (Right).}
\label{piceigen}
\end{figure}

The eigenvalues of the system \eqref{eq:semiDSM} can be obtained by the standard Fourier analysis 
 \cite{LiWang19}, given by  $\lambda_{j}=\frac{2D}{{\delta x}^2}(\cos(2\pi j\delta x)-1)-i \frac{a}{{\delta x}} \sin(2\pi j \delta x), j=1,2,...,N$, 
whose distributions and the corresponding numerical solutions are plotted in \cref{piceigen}.
We can see that the semi-discrete system has typical stiff characteristics and their eigenvalues 
have large non-zero imaginary parts, hence 
the numerical solutions exhibit oscillations and decays. To determine the decay rate of the numerical solutions, we can define the index function $p_{\alpha}$ as in \eqref{eq:indexp}, whose observed values are given in  \cref{table21} and \cref{table22}, from which we observe that the solution presents an algebraic decay rate and the index function $p_{\alpha}$ is in perfect agreement with our theoretical prediction.

\begin{table}\label{table21}
\caption{ Observed $p_\alpha$ computed by $\mathcal{L}1$ and F-BDF$1$(the data in the brackets) with $h=0.01$, $a=0.1$, $D=5$, $N=64$ and initial value $u_0=10\sin(4\pi x)$ for Example 5.2.} 
\begin{center}
\begin{tabular}{lll lll ll }
\hline $t_n$  & $\alpha=0.3$  & $\alpha=0.5$ & $\alpha=0.7$ & $\alpha=0.9$\\
\hline
$10   $ & 0.3003(0.3004)  & 0.5007(0.5009) & 0.7012(0.7014)  & 0.9016(0.9020)\\
$20   $ & 0.3001(0.3002)  &0.5004(0.5004)  & 0.7006(0.7007)  & 0.9008(0.9010)\\
$30  $ & 0.3001(0.3001)  &0.5002(0.5003) & 0.7004(0.7005)  & 0.9005(0.9007)\\
$40  $& 0.3000(0.3001)  & 0.5002(0.5002) & 0.7003(0.7004)  & 0.9004(0.9005)\\
$50 $& 0.3000(0.3000)  & 0.5001(0.5002) & 0.7003(0.7003)  & 0.9003(0.9004)\\
\hline
\end{tabular}
\end{center}
\end{table}

\begin{table}\label{table22}
\caption{ Observed $p_\alpha$ computed by F-BDF$2$ and F-Adams$2$ (the data in the brackets) with $h=0.01$, $a=0.1$, $D=5$, $N=64$ and initial value $u_0=10\sin(4\pi x)$ for Example 5.2.} 
\begin{center}
\begin{tabular}{lll lll ll }
\hline $t_n$  & $\alpha=0.3$  & $\alpha=0.5$ & $\alpha=0.7$ & $\alpha=0.9$\\
\hline
$10   $ & 0.3003(0.3003)  & 0.5007(0.5007) & 0.7012(0.7012)  & 0.9016(0.9016)\\
$20   $ & 0.3001(0.3001)  &0.5004(0.5004)  & 0.7006(0.7006)  & 0.9008(0.9008)\\
$30  $ & 0.3001(0.3001)  &0.5002(0.5002) & 0.7004(0.7004)  & 0.9005(0.9005)\\
$40  $& 0.3000(0.3000)  & 0.5002(0.5002) & 0.7003(0.7003)  & 0.9004(0.9004)\\
$50 $& 0.3000(0.3000)  & 0.5001(0.5001) & 0.7003(0.7003)  & 0.9003(0.9003)\\
\hline
\end{tabular}
\end{center}
\end{table}

\subsection{Fractional Lorenz controlled system} Consider the fractional Lorenz controlled system
\begin{equation}\label{eq:Fcs}
\begin{split}
_{~0}\mathcal{D}_{t}^{\alpha}y(t)&= Ay+f(y)+Bu,\q  u(t)=Ky(t),
\end{split}
\end{equation}
where the coefficient matrices $A$ and the nonlinear function $f$ are given by 
$$A=\begin{pmatrix}
-10 & 10 & 0 \\
28 & -1 & 0\\
0 & 0 & -\frac{8}{3} \\
\end{pmatrix},\quad
f(y)=\begin{pmatrix}
0\\
-y_{1}(t) y_3(t)\\
y_{1}(t) y_2(t)\\
\end{pmatrix}.
$$
When there is no control (i.e., $B\equiv0$), it is known the fractional Lorenz system has chaos solutions, 
similar to the classical Lorenz system, which are uniformly bounded but do not decay 
to some equilibrium points. It was shown in \cite{Cong19} that 
the trivial solution to this controlled system (with $B=(1,1,1)^{T}$ and $K=(0, -10, 0)$) 
is Mittag-Leffler stable for all $\alpha\in(0,1)$.

\begin{figure}[!ht]
\begin{center}
 \includegraphics[scale=0.42]{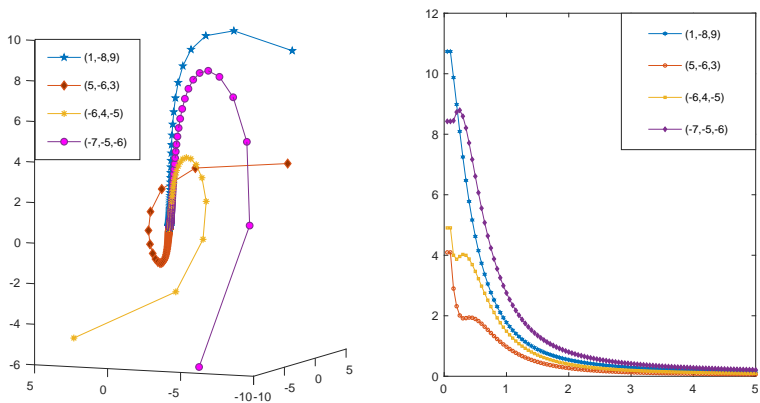}
\end{center}
\caption{The numerical solutions (left) and the norm $\|y_{n}\|$(right) for various initial values with parameters $\alpha=0.9, h=0.05$ computed by F-BDF$1$ for Example 5.3.}
\label{pic54}
\end{figure}

\begin{table}\label{table41}
\caption{Observed index $p_\alpha$ computed by $\mathcal{L}1$ and F-BDF$1$ (in the brackets) methods, 
with $h=0.1$ and the initial values $(1,-8,9)$ for Example 5.3.} 
\begin{center}
\begin{tabular}{lll lll ll }
\hline $t_n$  & $\alpha=0.3$  & $\alpha=0.5$ & $\alpha=0.7$ & $\alpha=0.9$\\
\hline
$20  $ & 0.2770(0.2771)  & 0.5026(0.5032) & 0.7335(0.7348)  & 0.9502(0.9525)\\
$40  $ & 0.2807(0.2808)  &0.5018(0.5021)  & 0.7199(0.7206)  & 0.9257(0.9267)\\
$60  $ & 0.2827(0.2828)  &0.5014(0.5016)  & 0.7147(0.7152)  & 0.9175(0.9182)\\
$80  $& 0.2840(0.2841)  & 0.5012(0.5013)  & 0.7119(0.7122)  & 0.9134(0.9139)\\
$100 $ & 0.2850(0.2850)  & 0.5010(0.5012)  & 0.7101(0.7104)  & 0.9109(0.9113)\\
\hline
\end{tabular}
\end{center}
\end{table}

\begin{table}\label{table42}
\caption{Observed index $p_\alpha$ computed by F-BDF$2$ and F-Adams$2$ (in the brackets) methods for $h=0.1$ and initial values $(1, -8, 9)$ for Example 5.3.} 
\begin{center}
\begin{tabular}{lll lll ll }
\hline $t_n$  & $\alpha=0.3$  & $\alpha=0.5$ & $\alpha=0.7$ & $\alpha=0.9$\\
\hline
$20   $& 0.2770(0.2770)  &0.5026(0.5026) & 0.7334(0.7334)  & 0.9499(0.9499)\\
$40   $& 0.2807(0.2807)  &0.5018(0.5018) & 0.7199(0.7199)  & 0.9256(0.9256)\\
$60  $ & 0.2827(0.2827)   &0.5014(0.5014) & 0.7147(0.7147)  & 0.9175(0.9175)\\
$80  $& 0.2840(0.2840)   &0.5012(0.5012)  & 0.7119(0.7119)  & 0.9133(0.9133)\\
$100  $& 0.2850(0.2850)   &0.5010(0.5011)  & 0.7101(0.7101)  & 0.9108(0.9108)\\
\hline
\end{tabular}
\end{center}
\end{table}

\begin{table}\label{table43}
\caption{The observed index functions $p_\alpha$ computed by $\alpha$-difference method with $h=0.1$ and initial values $(1,-8,9)$ for Example 5.3.} 
\begin{center}
\begin{tabular}{lll lll ll }
\hline $t_n$  & $\alpha=0.3$  & $\alpha=0.5$ & $\alpha=0.7$ & $\alpha=0.9$\\
\hline
$20  $ & 1.2508  &1.4989 & 1.7625  & 1.9993\\
$40  $ & 1.2579 &1.4995  & 1.7376  & 1.9502\\
$60  $ & 1.2621 &1.4996 & 1.7279  & 1.8578\\
$80  $& 1.2648   &1.4997 & 1.7226  & 1.8645\\
$100  $& 1.2664  &1.4998  & 1.7123  & 1.7928\\
\hline
\end{tabular}
\end{center}
\end{table}

The numerical solutions and the norm $\|y_{n}\|$ are plotted in \cref{pic54}, 
with different initial values. We can observe that under small perturbations, 
the numerical solutions always converge to the equilibrium point even they start with different initial values,
and the convergence rate exhibits a typical algebraic decay rate. The numerically observed indices $p_{\alpha}$ 
obtained by $\mathcal{L}1$ method, F-BDF$1$, F-BDF$2$,  F-Adams$2$ and $\alpha$-Difference method are presented in 
Tables \ref{table41}, \ref{table42} and \ref{table43},  respectively. 
The results in Tables \ref{table41} and  \ref{table42}  are almost consistent with our theoretical prediction with the decay rate  $\|y_n\|=O(t_n^{-\alpha})$.  When $\alpha$ is small, such as $\alpha=0.3$, the observed values are slightly lower than what we would expect. 
The reason is that when $\alpha$ is small, it usually takes much longer time for the system to decay into equilibrium.
The results in \cref{table43} have higher decay rate, which is about $\|y_n\|=O(t_n^{-1-\alpha})$.   
We emphasize that the solutions all maintain the typical polynomial decay rate, which is significantly different from the exponential decay of solutions to integer order equations.

\section{Concluding remarks}\label{sec:dis}
We have established the numerical Mittag-Leffler stability of the strongly A-stable F-LMMs and the $\mathcal{L}1$ method through the singularity analysis of the generating functions of numerical schemes for linear F-ODEs. This stability describes the optimal long-term algebraic decay rate of the numerical solutions, and it is shown both analytically and numerically that the algebraic decay rate of numerical solutions is exactly preserved 
as that of the continuous solutions.
For semi-linear F-ODEs with small perturbations, the numerical Mittag-Leffler stability near the equilibrium is also 
derived, by making use of some new and improved discrete resolvent operator estimates. 
These results have reconfirmed the slow diffusion of solutions to time fractional equations 
over a long period of time from a numerical point of view.

The analysis and theory developed in this work may help us better understand and analyze  
more complex time fractional PDEs. For example, it was proved in \cite{Zacher19} 
by the entropy method that the solutions to the time fractional Fokker-Planck equations
converge to an equilibrium state in $L^1$-norm with an algebraic decay rate;
it was observed through numerical experiments \cite{Tang19} that 
the numerical solutions to the time fractional phase field model exhibits 
an algebraic decay rate and slow energy dissipation and the solutions of 
time fractional Allan-Chan equations decay as $O(t^{-\frac{\alpha}{3}})$ in $L^2$-norm. 
But there are still no rigorous numerical analysis for these models and asymptotic behaviors of their solutions.

\appendix
\section{Fractional calculus, Mittag-Leffler function and Prabhakar function}
\label{sec:appendixA}
The Caputo fractional derivative of order $\alpha\in(0, 1)$ is given by \cite{Pod}
\begin{equation*} 
\begin{split}
\mathcal{D}_{t}^{\alpha}y(t):= \mathcal{I}_{t}^{1-\alpha}y'(t)=(k_{1-\alpha}*y')(t)=\frac{1}{\Gamma(1-\alpha)}\int_0^t \frac{y'(s)}{(t-s)^{\alpha}}ds,  t>0,
\end{split}
\end{equation*}
where $\mathcal{I}_{t}^{\alpha}y(t)=(k_{\alpha}*y)(t)=\frac{1}{\Gamma(\alpha)}\int_0^t \frac{y(s)}{(t-s)^{1-\alpha}}ds$ denotes the Riemann-Liouville integral and the stand kernel $k_{\alpha}(t)=\frac{t^{\alpha-1}}{\Gamma(\alpha)}$.  
We recall the Mittag-Leffler functions $E_{\alpha}(z)$
and $E_{\alpha,\beta}(z)$:
$E_{\alpha}(z)=\sum_{k=0}^{\infty}\frac{z^{k}}{\Gamma(\alpha
k+1)},\, E_{\alpha,\beta}(z)=\sum_{k=0}^{\infty}\frac{z^{k}}{\Gamma(\alpha
k+\beta)}$, where $\alpha,\beta>0$ and $z\in \mathbb{C}$, 
which can be seen as the fractional generalization of exponential functions and occur naturally 
in fractional calculus; see more details in \cite{Pod}.
For $\alpha\in (0,1)$, these two functions have the asymptotic expansion 
\begin{equation}\label{eq:MLexpan}
\begin{split}
E_{\alpha,\beta}(z)=-\sum_{k=1}^{N}\frac{1}{\Gamma(\beta-k\alpha)}\frac{1}{z^{k}}+O\left(\frac{1}{|z|^{N+1}}\right), ~~N\in\mathbb{N}^{+},\, |z|\to\infty, \, \theta\leq \arg(z)\leq\pi
\end{split}
\end{equation}
where $\theta\in(\frac{\alpha\pi}{2}, \pi\alpha)$.
According to the expansion \eqref{eq:MLexpan} 
one can prove that \cite{Cong18,Cong19}
\begin{equation}\label{eq:MLexpan3}
\begin{split}
\left| E_{\alpha}(\lambda t^\alpha)\right| \leq \frac{C_1(\alpha,\lambda)}{t^{\alpha}}, \qquad 
\left| E_{\alpha,\alpha}(\lambda t^\alpha)\right| \leq \frac{C_2(\alpha,\lambda)}{t^{2\alpha}}, ~~\forall t\geq t_0>0,
\end{split}
\end{equation}
where $\lambda\in \Lambda_{\alpha}^{s}$ and $C_1(\alpha,\lambda), C_2(\alpha,\lambda)$ are real positive constants which are independent of $t$.

The Prabhakar function is the three-parameter generalization of Mittag-Leffler function defined as
\begin{equation*}
\begin{split}
E_{\alpha,\beta}^{\gamma}(z)=\sum_{k=0}^{\infty}\frac{(\gamma)_{k}z^{k}}{k! \Gamma(\alpha
k+\beta)}, \quad z\in \mathbb{C},
\end{split}
\end{equation*}
where $(\gamma)_{k}=\Gamma(\gamma+k)/\Gamma(k)$ is the Pochhammer symbol. 
It is enough for our purpose to restrict the parameters $\alpha,\beta, \gamma\in\mathbb{R} ~\text{and}~\alpha>0$.
In this case, $E_{\alpha,\beta}^{\gamma}(z)$ is an entire function of order $\rho= 1/\alpha$.
The key historical events and modern development and the main properties of Prabhakar function can be found in the most recent review paper
\cite{Garrappa20}. We have $E_{\alpha,\beta}^{\gamma}(z)=E_{\alpha,\beta}(z)$ 
for $\gamma=1$, which recovers the Mittag-Leffler function.
Generally, we have the reduction formula 
\begin{equation}\label{eq:Prafun}
\begin{split}
E_{\alpha,\beta}^{\gamma+1}(z)=\frac{ E_{\alpha,\beta-1}^{\gamma}(z)+(1-\beta+\alpha\gamma) E_{\alpha,\beta}^{\gamma}(z) }{\alpha\gamma}.
\end{split}
\end{equation}

\section{ Estimate for $D_n$ in \eqref{eq:dcvf3} based on standard resolvent \eqref{eq:resestimate}}
\label{sec:appendixB}

The integral path in \eqref{eq:dcvf3} is firstly deformed to $\Gamma_{(r, \theta)}$ and then split into two parts, 
$\Gamma_{1}$ and $\Gamma_{2}$, where $\Gamma_{1}$ is the circle around the origin and $\Gamma_{2}$ is two line segments. Then the integral is naturally divided into $D_{n} =D^{1}_{n}+ D^{2}_{n}$,  where $D_{n}^{j}$ $(j=1, 2)$ means the integral carried over $z\in\Gamma_{j}$. One of the main reasons for choosing this contour is that the exponential function $e^{z t_{n}}$ will decay as $t_{n}$ increases for $z\in\Gamma_{2}$. 
In order to estimate $D^{j}_{n}$, we need to estimate three terms in the integral both for $z\in\Gamma_{1}$ and  $z\in\Gamma_{2}$, namely, 
$|e^{z t_{n}}|, \|(h^{-\alpha}  F_{\mu} (e^{-zh}) I- A)^{-1}\|$ and $|dz|$, where 
$F_{\mu} (z)=(1-z)^\alpha$ for F-BDF$1$.

We first consider $D^{1}_{n}$ for $z\in \Gamma_{1}$. 
Let $z=re^{i\theta}=\frac{1}{t_{n}}e^{i\theta}$, where $-\phi\leq \theta \leq \phi$ and  $\phi\in(\pi/2, \pi)$. 
 By direct calculation we can get that 
 $|e^{z t_{n}}|=|e^{\frac{1}{t_{n}}(\cos(\theta)+i\sin(\theta)) t_{n}}|= |e^{\cos(\theta)} |\leq e$ 
 and $\int_{\Gamma_{1}}|dz|=\ell(\Gamma_{1}) \leq 2\pi r\leq \frac{C}{t_{n}}$, where $\ell(\Gamma^{1})$ is the length of $\Gamma_{1}$.
 According to \eqref{eq:resestimate}, we have $\| (h^{-\alpha}  F_{\mu} (e^{-zh}) I- A)^{-1}\|\leq C |\frac {h^{\alpha}} {(1-e^{-zh})^{\alpha}}|$.
On the other hand, there exist constants $c_1, c_2>0$ such that $c_{1}|zh|\leq |1-e^{-zh}|\leq c_{2}|zh|$ for $z\in \Gamma_{1}$ \cite{Jin17}, which yields that  
$|\frac {h^{\alpha}} {(1-e^{-zh})^{\alpha}}| \leq \frac{C}{|z|^{\alpha}}=\frac{C}{t_{n}^{\alpha}}.$
Combining the above bounds lead to
\begin{equation} \label{eq:estimate4}
\begin{split}
\|D^{1}_{n} \| \leq  \frac{e}{2\pi } \cdot \frac{C}{t_{n}^{\alpha}}\cdot \frac{C}{t_{n}} =O(t_{n}^{-\alpha-1}).
\end{split}
\end{equation}

We now estimate $D^{2}_{n}$. Let $z=re^{i\phi}$, where ${1}/{t_{n}}\leq r\leq \pi/(h\sin(\phi))$ and $\phi\in(\pi/2, \pi)$.
It follows by direct calculation that 
$|e^{z t_{n}}|=|e^{r(\cos(\phi)+i\sin(\phi)) t_{n}}| =|e^{r\cos(\phi) t_{n}}|= e^{r\cos(\phi) t_{n}}$.
Similarly to the above estimation of $D^{1}_{n}$, there exist constants $c_1, c_2>0$ such that $c_{1}|zh|\leq |1-e^{-zh}|\leq c_{2}|zh|$ for $z\in \Gamma_{2}$, which yields that  
$|\frac {h^{\alpha}} {(1-e^{-zh})^{\alpha}}| \leq \frac{C}{|z|^{\alpha}}=\frac{C}{r^{\alpha}}.$
Combining the above bounds and the simple transformation $s=rt_{n}$, 
we can derive
\begin{equation} \label{eq:estimate32}
\begin{split}
\|D^{2}_{n} \| 
&\leq \frac{1}{2\pi } \int_{\frac{1}{t_{n}}}^{\frac{\pi}{h\sin(\phi)}} e^{r\cos(\phi) t_{n}} \cdot \frac{C}{r^{\alpha}}  dr
= \frac{C}{2\pi }   \int_{\frac{t_{n}}{t_{n}}}^{\frac{\pi}{h\sin(\phi)} t_{n}} e^{s\cos(\phi)} \cdot  \frac{t_{n}^\alpha }{s^{\alpha}} \cdot \frac{1}{t_n}ds \\
&=  \frac{C}{2\pi }   \cdot \frac{1}{t_{n}^{1-\alpha}}  \int_{1}^{\frac{\pi}{h\sin(\phi)} t_{n}} e^{s\cos(\phi)} s^{-\alpha} ds=O(t_{n}^{\alpha-1}). 
\end{split}
\end{equation}
Therefor we conclude that $\|D_{n}\|\leq \|D^{1}_{n} \| +\|D^{2}_{n} \| =O(t_{n}^{\alpha-1})$ as $t_{n}\to \infty.$

From the above analysis, we see that $\|D^{1}_{n} \| =O(t_{n}^{-\alpha-1})$.  However, for $D^{2}_{n}$ with $z\in \Gamma_{2}$, this method only gives that $\|D^{2}_{n} \|=O(t_{n}^{\alpha-1})$, hence resulting in 
the estimate $\|D_{n}\|=O(t_{n}^{\alpha-1})$ or a reduction of the decay rate.
In this case, the sum of the series  $\sum_{n=1}^{\infty}\|D_{n}\|$ diverges.
This is not satisfactory and can not be used to establish the long-time decay of numerical solutions.

\section*{Acknowledgements}
The authors are grateful to Professor Lei Li (Shanghai Jiao Tong University) for many constructive discussions.

\end{document}